\documentclass[12pt]{article}
\usepackage{amsmath, amsthm}
\usepackage{amscd}
\usepackage{amssymb}
\usepackage{array}

\newtheorem{thm}{Theorem}[section]
\newtheorem{theorem}[thm]{Theorem}

\newtheorem{lemma}[thm]{Lemma}
\newtheorem{proposition}[thm]{Proposition}

\theoremstyle{definition}
\newtheorem{definition}[thm]{Definition}
\newtheorem{example}[thm]{Example}

\newtheorem{observation}[thm]{Observation}

\newtheorem{remark}[thm]{Remark}

\begin{document}

\newcommand{\id}{\relax{\rm 1\kern-.28em 1}}
\newcommand{\R}{\mathbb{R}}
\newcommand{\C}{\mathbb{C}}
\newcommand{\Z}{\mathbb{Z}}
\newcommand{\g}{\mathfrak{G}}
\newcommand{\fp}{\mathfrak p}
\newcommand{\fn}{\mathfrak n}
\newcommand{\fh}{\mathfrak h}
\newcommand{\fs}{\mathfrak s}
\newcommand{\fc}{\mathfrak c}
\newcommand{\fk}{\mathfrak k}
\newcommand{\fl}{\mathfrak l}
\newcommand{\frs}{\mathfrak s}
\newcommand{\e}{\epsilon}

\newcommand{\bp}{\mathbf{p}}
\newcommand{\bmax}{\mathbf{m}}
\newcommand{\bT}{\mathbf{T}}
\newcommand{\bU}{\mathbf{U}}
\newcommand{\bP}{\mathbf{P}}
\newcommand{\bA}{\mathbf{A}}
\newcommand{\bm}{\mathbf{m}}
\newcommand{\bIP}{\mathbf{I_P}}

\newcommand{\cU}{\mathcal{U}}
\newcommand{\cA}{\mathcal{A}}
\newcommand{\cB}{\mathcal{B}}
\newcommand{\cT}{\mathcal{T}}
\newcommand{\cI}{\mathcal{I}}
\newcommand{\cO}{\mathcal{O}}
\newcommand{\oo}{\mathcal{O}}
\newcommand{\cG}{\mathcal{G}}
\newcommand{\cJ}{\mathcal{J}}
\newcommand{\cF}{\mathcal{F}}
\newcommand{\cK}{\mathcal{K}}
\newcommand{\dd}{\mathcal{D}}
\newcommand{\E}{\mathcal{E}}
\newcommand{\cH}{\mathcal{H}}
\newcommand{\cM}{\mathcal{M}}
\newcommand{\cPO}{\mathcal{PO}}
\newcommand{\cl}{\ell}
\newcommand{\cFG}{\mathcal{F}_{\mathrm{G}}}
\newcommand{\cHG}{\mathcal{H}_{\mathrm{G}}}
\newcommand{\Gal}{G_{\mathrm{al}}}
\newcommand{\rLie}{{\mathrm{Lie}}}

\newcommand{\rGL}{\mathrm{GL}}
\newcommand{\rSU}{\mathrm{SU}}
\newcommand{\rSL}{\mathrm{SL}}
\newcommand{\rSO}{\mathrm{SO}}
\newcommand{\rOSp}{\mathrm{OSp}}
\newcommand{\rsl}{\mathrm{sl}}
\newcommand{\rM}{\mathrm{M}}
\newcommand{\M}{\mathrm{M}}
\newcommand{\End}{\mathrm{End}}
\newcommand{\Hom}{\mathrm{Hom}}
\newcommand{\diag}{\mathrm{diag}}
\newcommand{\rspan}{\mathrm{span}}
\newcommand{\rank}{\mathrm{rank}}
\newcommand{\Gr}{\mathrm{Gr}}
\newcommand{\ber}{\mathrm{Ber}}

\newcommand{\fsl}{\mathfrak{sl}}
\newcommand{\fg}{\mathfrak{g}}
\newcommand{\ff}{\mathfrak{f}}
\newcommand{\fgl}{\mathfrak{gl}}
\newcommand{\fosp}{\mathfrak{osp}}
\newcommand{\fm}{\mathfrak{m}}
\newcommand{\fso}{\mathfrak{so}}
\newcommand{\fsu}{\mathfrak{su}}

\newcommand{\str}{\mathrm{str}}
\newcommand{\Sym}{\mathrm{Sym}}
\newcommand{\tr}{\mathrm{tr}}
\newcommand{\defi}{\mathrm{def}}
\newcommand{\Ber}{\mathrm{Ber}}
\newcommand{\spec}{\mathrm{Spec}}
\newcommand{\sschemes}{\mathrm{(sschemes)}}
\newcommand{\sschemeaff}{\mathrm{ {( {sschemes}_{\mathrm{aff}} )} }}
\newcommand{\set}{\mathrm{(set)}}
\newcommand{\Top}{\mathrm{Top}}
\newcommand{\sarf}{ \mathrm{ {( {salg}_{rf} )} }}
\newcommand{\smfld}{\mathrm{ {(smfld)} }}
\newcommand{\odd}{\mathrm{odd}}
\newcommand{\alg}{\mathrm{(alg)}}
\newcommand{\sa}{\mathrm{(salg)}}
\newcommand{\SA}{\mathrm{(salg)}}
\newcommand{\salg}{\mathrm{(salg)}}
\newcommand{\varaff}{ \mathrm{ {( {var}_{\mathrm{aff}} )} } }
\newcommand{\svaraff}{\mathrm{ {( {svar}_{\mathrm{aff}} )}  }}
\newcommand{\ad}{\mathrm{ad}}
\newcommand{\Ad}{\mathrm{Ad}}
\newcommand{\pol}{\mathrm{Pol}}
\newcommand{\Lie}{\mathrm{Lie}}
\newcommand{\Proj}{\mathrm{Proj}}
\newcommand{\uspec}{\underline{\mathrm{Spec}}}
\newcommand{\uproj}{\mathrm{\underline{Proj}}}
\newcommand{\rGr}{\mathrm{Gr}}

\newcommand{\sym}{\cong}

\newcommand{\al}{\alpha}
\newcommand{\be}{\beta}
\newcommand{\lam}{\lambda}
\newcommand{\de}{\delta}
\newcommand{\ga}{\gamma}
\newcommand{\lra}{\longrightarrow}
\newcommand{\ra}{\rightarrow}

\newcommand{\bM}{{\bar{M}}}
\newcommand{\bz}{\bar z}
\newcommand{\btheta}{\bar \theta}


\rightline{IFIC/06-28}

\rightline{FTUV-05/0928}

\smallskip

    \centerline{\LARGE \bf  The Minkowski and }

    \smallskip

     \centerline{\LARGE \bf conformal superspaces}

\vskip 1cm

\centerline{ R. Fioresi \footnote{Investigation supported by the
University of Bologna, funds for selected research topics.}}

\smallskip

\centerline{\it Dipartimento di Matematica, Universit\`a di Bologna
}
 \centerline{\it Piazza di Porta S. Donato, 5.}
 \centerline{\it 40126 Bologna. Italy.}
\centerline{{\footnotesize e-mail: fioresi@dm.UniBo.it}}

\medskip

\centerline{ M. A. Lled\'o }

\smallskip

 \centerline{\it  Departament de F\'{\i}sica Te\`orica,
Universitat de Val\`encia and IFIC}
 \centerline{\small\it C/Dr.
Moliner, 50, E-46100 Burjassot (Val\`encia), Spain.}
 \centerline{{\footnotesize e-mail: Maria.Lledo@ific.uv.es}}

\medskip

\centerline{and V. S. Varadarajan}

 \smallskip

\centerline{\it Department of Mathematics, UCLA.} \centerline{\it
Los Angeles, CA, 90095-1555, USA} \centerline{{\footnotesize e-mail:
vsv@math.ucla.edu}}

 \vskip 2cm

\begin{abstract}
  We define complex Minkowski superspace in 4 dimensions as the big cell
inside a complex flag supermanifold. The complex conformal
supergroup acts naturally on this super flag, allowing us to
interpret it as the conformal compactification of complex Minkowski
superspace. We then consider real
 Minkowski superspace as a suitable real form of the complex version. Our
methods are group theoretic, based on the real  conformal supergroup
and its
 Lie superalgebra.
\end{abstract}

\vfill\eject

\tableofcontents

\vfill\eject

\section{Introduction}

The superfield approach to supersymmetric physical theories is
based on Minkowski superspace time. In the classical world the
Poincar\'{e} group, which acts on Minkowski space, is a symmetry of
any physically sensible theory. The Poincar\'{e} group is a
non-semisimple group. The smallest simple group containing the
Poincar\'{e} group in dimension 4 is the conformal group $\rSO(4,2)$,
or its double covering $\rSU(2,2)$. Conformal symmetry is usually
broken in nature. But fields of zero mass are conformally
invariant, so there may be regimes where this symmetry can be
restored. Minkowski space time is not enough to support conformal
symmetry, and it needs to be {\it compactified} by adding points
at infinity to obtain a space with an action of the conformal
group. It turns out that the compactified Minkowski space is in
fact a real form of the complex Grassmannian of two planes in a
complex four dimensional vector space (see Ref. \cite{ma1} for a
careful treatment of this picture.)

The complex Grassmannian is a homogeneous space for the complex group $\rSL(4,\C)$ and so
has a cellular decomposition derived from the Bruhat decomposition of the complex group.
Complex Minkowski space time is thus the {\it big cell} inside the Grassmannian. One can go
from conformal symmetry to Poincar\'{e} symmetry and
 {\it vice versa} by deleting or adding the points at infinity. In a way,
 the conformal description may be useful even when the conformal
 group is not a symmetry. This was one of the justifications behind
 the seminal paper by Penrose about what he called  the  {\it twistor}
 approach to Minkowski space time \cite{pe}.

In the super world this picture can be generalized completely.
Supertwistors were introduced in Ref. \cite{fe}. The purpose of this
paper is to give a mathematically careful description of this
generalization. We concentrate on Minkowski space time of dimension
4 and signature (1,3)\footnote{This is the signature of the flat
metric in Minkowski space.} and consider its extension to a
superspace of odd dimension 4, which means, in terms known to
physicist, that we have $N=1$ supersymmetry carried by a Majorana
spinor of four components.

The case of even dimension four is singled out for physical reasons.
Most of the treatment can be appropriately generalized to higher
dimensions (and other signatures). Nevertheless the  occurrence of
Grassmannian and flag manifolds is peculiar to dimension four. The
extension to higher supersymmetries does not present a problem, and
we have chosen $N=1$ for concreteness.

The discovery that superconformal spacetime is a {\it flag} super
manifold, rather than a Grassmannian supermanifold, as one may
naively be led to believe from what happens in the classical
case\footnote{A Grassmannian supermanifold would be useful in
describing {\it chiral fields} \cite{fwz}, which are not real.}, is
due to Manin \cite{ma1} and Kotrla and Niederle \cite{kn}. Working
inside a flag super manifold is necessary to obtain the correct real
form. This marks a deep difference between the classical and the
super versions. We refer also to Ref. \cite{hh} for a quite complete
review on superspaces associated with four dimensional Minkowski
spacetime.

It is our goal to give a group theoretical view of this
construction. While it is standard in the classical case, the
supergroup version presents many peculiarities that have to be
worked out. We construct the super flag explicitly as a quotient
space of the superconformal group, keeping always the action of the
supergroup explicit.  This leads us to a general discussion of
homogeneous spaces of supergroups which had not been rigourously
treated in the literature before. We obtain the homogeneous space
$G/H$, with $G\subset H$, both supergroups, as the unique
supermanifold satisfying some natural properties.

The definition of Minkowski and conformal superspaces is not
arbitrary, since they are determined by physical requirements. In
fact, we provide in precise form the requirements that determine
uniquely the Poincar\'{e} and translation superalgebras inside the
conformal superalgebra, both over the complex and the real fields.
The correct real form is obtained as the fixed points of a
conjugation that we give explicitly. This solves the problem at the
infinitesimal level, since the dual of the translation superalgebra
can be identified with the infinitesimal Minkowski superspace.

The passage to the global theory is done by constructing the flag
supermanifold as an homogeneous superspace, reducing infinitesimally
to the Lie algebra decomposition found before. The Minkowski
superspace is then the big cell inside the flag supermanifold. The
real form is obtained by lifting the conjugation found in the Lie
superalgebra to the supergroup and then reducing it to the
homogeneous space. We also obtain the real form of the Minkowski
superspace from the real form of the flag suepermanifold. This
provides a significantly different and more natural description than
the one appearing in the above mentioned papers \cite{ma1,kn}.

We want the paper to be as self contained as possible, so we include
a review of the classical case in Section \ref{minkowski} and we
explain the language of supergeometry in some detail in the Appendix
\ref{supergeometry}. Also, we have tried to work out every result
explicitly in coordinates, as to make contact with the physics
notation.

\medskip

The paper is organized as follows.

In Section \ref{minkowski} we review the construction of classical
Minkowski space time in four dimensions. We start with the complex
version, viewed as the big cell inside the Grassmannian of two
planes in $\C^4$. The complex conformal group is $\rSL(4,\C)$, which
acts naturally on the Grassmannian. One can see the complexified
Poincar\'{e} group (including dilations) as the subgroup of $\rSL(4,\C)$
leaving the big cell invariant. Then the Grassmannian is the
conformal compactification of the complex Minkowski space. We then
describe the conjugation that leads to the correct real form of the complex
Minkowski space time.

In Section \ref{homogeneous} we take the classical theory of
homogeneous spaces and extend it to the super setting. We make a
precise definition of what a homogeneous space of a supergroup is
and, using the results of Appendices \ref{supergeometry1} and
\ref{supergeometry2}, we give to it the structure of supermanifold,
together with an action of the supergroup on it. We also find its
functor of points.

In Section \ref{spst} we consider the problem at the infinitesimal
level. The  Poincar\'{e} Lie superalgebra can be defined through a set
of natural conditions. Then we show that there exists a unique
subalgebra of the conformal superalgebra satisfying these
conditions. In fact, the conformal superalgebra splits as the direct
sum (as vector spaces) of the Poincar\'{e} superalgebra plus another
subalgebra, the translation superalgebra of Minkowski space time, as
it arises in physics (from where it takes the name).

The infinitesimal treatment is crucial to find the adequate real
form, which is not completely straightforward in the super case.
When passing to the global theory in Section \ref{global} we will
see that, in order to obtain the desired real form, it is not
enough to consider a super Grassmannian, but we need a super flag
manifold \cite{ma1, kn}. In fact, one could generalize the
classical construction naively in the complex case, but the
reality properties of the spinor representations prevent it from
having a real structure. These peculiarities depend heavily on the
dimension and signature of the Minkowski space (see for example
Ref. \cite{dflv} for an account of super Poincar\'{e} and conformal
superalgebras in arbitrary dimension and signature).

We can lift the conjugation found at infinitesimal level to the
supergroup. In particular, it preserves the super Poincar\'{e} group,
and it descends to the superflag seen as a homogeneous space. It
also preserves the big cell thus defining the standard Minkowski
superspace.

In Appendices \ref{supergeometry1} and \ref{supergeometry2} we
recall some facts on superspaces and supermanifolds which will be
needed through the discussion. We emphasize the functorial approach:
a supermanifold can be given in terms of a representable functor,
its {\it functor of points}. In particular, we provide a
representability criterion which allows us to determine if a functor
is the functor of points of a supermanifold. This result is new.  We
illustrate the general theory with two examples that are essential
for us: the  Grassmannian and the  flag supermanifolds.

In Appendix \ref{supergroups} we give a brief account on super Lie
groups. Finally, given the importance of real forms of
supermanifolds in this work, we devote Appendix \ref{realform} to
defining them carefully.

\section{The complex Minkowski space time and the Grassmannian G(2,4)}
\label{minkowski}

This section is a summary of the facts that we will need from Ref.
\cite{va}, Section 3.3. All the details and proofs of our
statements, and many others, can be found there.

The basic idea is to describe complex Minkowski space time as a big
cell of the Grassmannian $G(2,4)$ by finding the action of the
Poincar\'e group as the subgroup of $\rGL(4,\C)$ which leaves the
big cell invariant.

We consider the Grassmannian $G(2,4)$, the set of two dimensional
subspaces or planes in $\C^4$. A plane can be  given by two linearly
independent vectors and we will denote it as
$\pi=(a,b)=\rspan\{a,b\}$, where $a$ and $b$ are column vectors.
Obviously if $\rspan\{a,b\}=\rspan\{a',b'\}$ they define the same
point of the Grassmannian. There is a transitive action of
$\rGL(4,\C)$ on $G(2,4)$,
$$g\in \rGL(4,\C),\qquad g\pi=(ga,gb).$$

Let $\{e_1,e_2,e_3, e_4\}$ be a basis of $\C^4$ and consider the
plane $\pi_0$ spanned by the vectors $e_1$ and $e_2$. The isotropy
group of $\pi_0$ is
$$P_0=\left\{\begin{pmatrix}L&M \\
0&R\end{pmatrix}\in \rGL(4,\C)\right\},$$ with
$L, M, R$ being $2\times 2$ matrices, $L$ and $R$ invertible. We
have then that $G(2,4)\approx\rGL(4,\C)/P_0$.

It is possible to use $\rSL(4,\C)$ instead, since its action is also
transitive.

\smallskip

Since the  vectors $a$ and $b$ spanning a plane are linearly
independent, the $4\times 2$ matrix
$$\pi=(a,b)=\begin{pmatrix}a_1&b_1\\a_2&b_2\\
a_3&b_3\\a_4&b_4\end{pmatrix}$$ has rank 2, and at least one of the
6 independent minors $y_{ij}=a_ib_j-b_ia_j$, $i \neq j$
is different
from zero. The set \begin{equation}U_
{ij}=\left\{(a,b)\;/\;y_{ij}\neq
0\right\}\label{openset}\end{equation} is an open set of $G(2,4)$.
It is easy to see that a plane in $U_{ij}$ can always be represented
by
\begin{eqnarray*}&a=e_i+\alpha e_k+\gamma e_l,\qquad b=e_j+\beta e_k+\delta
e_l,\\ &\hbox{with } k\neq l \hbox{ and } k,l\neq
i,j.\end{eqnarray*} The six open sets
$\left\{U_{ij}\right\}_{i<j}$ are a cover of $G(2,4)$, and in each
open set the numbers $(\alpha,\beta,\gamma,\delta)$ are
coordinates, so $U_{ij}\approx\C^4$. It will be convenient to
organize the coordinates as a $2\times 2$ matrix
\begin{equation}\label{coordinates}A=\begin{pmatrix}\alpha&\beta\\\gamma&\delta\end{pmatrix},
\qquad \pi=\begin{pmatrix}\id\\A\end{pmatrix},\end{equation} so
$U_{ij}\approx M_2(\C)$.

\subsection{The Pl\"ucker embedding and the Klein quadric}

We are going to see how $G(2,4)$ is embedded in the projective space
$\C P^5$. We consider the vector space
$E=\Lambda^2(\C^4)\thickapprox \C^6$, with basis $\{e_i\wedge
e_j\}_{i<j}$, $i,j=1,\dots , 4$.  For a plane $\pi=(a,b)$ we have
$$a\wedge b=\sum_{i<j}y_{ij}e_i\wedge e_j.$$
A change of basis $(a',b')=(a,b)u$, where $u$ is a complex $2\times
2$ non singular matrix, produces a change
$$a'\wedge b'=\det (u) a\wedge b,$$
so the image $\pi\rightarrow [a\wedge b]$ is well defined into the
projective space $P(E)\thickapprox \C P^5$. It is called the
Pl\"ucker map and it is not hard to see that it is an embedding. The
quantities $y_{ij}$ are the {\it homogeneous Pl\"ucker coordinates},
and the image of the Pl\"ucker map can be identified with the
solution of $p\wedge p=0$ for $p\in E$. In coordinates, this reads
as
\begin{equation}
y_{12} y_{34}+y_{23} y_{14}+y_{31}
y_{24}=0, \qquad y_{ij}=-y_{ji}. \label{klein}
\end{equation}
Equation (\ref{klein}) is a quadric in $\C P^5$, the Klein quadric
$K$, and we have just seen that $K=G(2,4)$.

Let $Q$ be the quadratic form in $E$ defining the Klein quadric,
$$Q(y)=y_{12} y_{34}+y_{23} y_{14}+y_{31}
y_{24}.$$ For any $g\in \rGL(4,\C)$ we have that $Q(gy)=\det(g)
Q(y)$, so if $g\in \rSL(4,\C)$ then  it preserves the quadratic form
$Q$. Then, the action of $\rSL(4,\C)$ on
$\C^6\approx\lambda^2(\C^4)$ defines a map $\rSL(4,\C)\rightarrow
\rSO(6,\C)$. Its kernel
  is $\pm 1$. Since
both groups have the same dimension this means that $\rSL(4,\C)$ is
 the spin group of $\rSO(6,\C)$.

For concreteness, let  us consider $U_{12}$. A plane $\pi\in U_{12}$
can always be represented as
$$\pi=\begin{pmatrix}\id\\A\end{pmatrix},
\qquad A=\begin{pmatrix}\alpha&\beta\\\gamma&\delta\end{pmatrix}.$$
In particular, $\pi_0=\rspan\{e_1, e_2\}$ has $A=0$. The Pl\"ucker
coordinates are
\begin{equation}(y_{12}, y_{23},y_{31},y_{14},y_{24},y_{34})=
(1,-\alpha,-\beta,\delta,-\gamma,\alpha\delta-\beta\gamma).
\label{pluckerbigcell}
\end{equation}

Let us denote by $\pi_\infty$ the plane spanned by $e_3$ and $e_4$,
so $\C^4=\pi_0\oplus \pi_\infty$. We say that a plane is finite if
it has no intersection with $\pi_\infty$. The set $U_{12}$ can also
be characterized as the set of finite planes, and it is a dense open
set in $G(2,4)\simeq K$. It is called the  {\it big cell} of
$G(2,4)$, and we will denote it by $K^\times$.

We describe now the complement of $U_{12}$ in $G(2,4)$. A plane that
has non-empty intersection with $\pi_\infty$ can always be
represented by two independent vectors $a\in \pi_\infty$ and $b$,
$$a=a_3e_3+a_4e_4,\qquad b=b_1e_1+b_2e_2+b_3e_3+b_4e_4,$$ so
$y_{12}=0$. Substituting in (\ref{klein}) we obtain
\begin{equation}y_{23} y_{14}+y_{31} y_{24}=0.
\label{cone}\end{equation} This represents the
closure of a cone in $\C^4$. Let us see this:

\begin{enumerate}
 \item If $y_{34}\neq 0$ we may choose $y_{34}=1$ and the
remaining coordinates may take any value, except for the
constraint (\ref{cone}). This is an affine cone in $\C^4$.
\item If $y_{34}= 0$, then one of the remaining $y$'s must be
different from zero. The solutions of (\ref{cone}) are taken up to
multiplication by a complex number, so we have a quadric in $\C
P^3$. With these points we obtain the closure of the cone in 1.
 \end{enumerate}

We are going to see that the complex Poincar\'e group sits inside
$\rGL(4,\C)$ as the subgroup $P$ that leaves $U_{12}$ invariant. So
we will have a natural identification of $U_{12}$ with the complex
Minkowski space time.  Consequently, $G(2,4)$ is a compactification
of Minkowski space time by the addition of the cone described above.

In fact, one can see that $P=P_\infty$, the group that leaves
$\pi_\infty$ invariant. In the representation that we have used,
$$\pi_\infty=\begin{pmatrix}0\\\id\end{pmatrix},$$
and $P_\infty$   can be written as
$$P_\infty=\left\{\begin{pmatrix}L&0\\NL&R\end{pmatrix},\;\; L, R
\;\;\mathrm{invertible}\right\}.$$ The bottom left entry is an
arbitrary $2\times 2$ matrix. It can be written  as
 $NL$, with $N$ arbitrary, since $L$ is invertible.

The action on $U_{12}$ is
\begin{equation}A\mapsto N+RAL^{-1},\label{poincare}\end{equation} so $P$ has the structure of semidirect
product $P=M_2\ltimes H$, where $M_2$, the set of $2\times 2$
matrices $N$, are the translations and
$$H=\left\{\begin{pmatrix}L&0\\0&R\end{pmatrix},\; \; L,R\in
\rGL(2,\C),\;\; \det\!L\cdot \det\!R=1\right\}.$$ The subgroup $H$
is the direct product  $\rSL(2,\C)\times \rSL(2,\C)\times
\C^\times$. But $\rSL(2,\C)\times \rSL(2,\C)$ is the spin group of
$\rSO(4,\C)$, the complexified Lorentz group, and $\C^\times$ acts
as a dilation.

$G(2,4)$ is then a compactification of the Minkowski space time that is
equivariant under the action of the conformal group.

\subsection{The real form}

The real Minkowski space time is obtained by taking the set of
fixed points under a certain conjugation of the complex Minkowski
space time. This can be done compatibly with our construction. We
first consider the following  conjugation of $E=\Lambda^2(\C^4)$:
$$\theta: (y_{12}, y_{23},y_{31},y_{14},y_{24},y_{34})\;\mapsto
\;(\bar y_{12}, \bar y_{23}, \bar y_{24}, \bar y_{14}, \bar
y_{31}, \bar y_{34}). $$ The conjugation $\theta$ is also well
defined when passing to the projective space  $P(E)$. The
quadratic form $Q$ defining the  Klein quadric satisfies
$$\overline{Q(y)}=Q\left (\theta (y)\right).$$
The set of fixed points of $\theta$,
\begin{gather*}E^\theta=\{y \in E\;\;|\;\;
\theta(y)=y\}=\{y\in \Lambda^2(\C^4)\;/\; \\y_{12}, y_{23},
y_{34}, y_{14}\,\hbox{are real and }\, y_{31}=\bar
y_{24}\}\end{gather*} is a real vector space. We denote by $Q_R$
the quadratic form $Q$ restricted to $E^\theta$. The form $Q_R$ is
a real quadratic form with signature (4,2):
$$Q_R(y)= y_{12} y_{34}+y_{23}
y_{14}+y_{31} \bar y_{31}.$$

Let $K=\{[y]\in P(E)\;|\; Q(y)=0\}$. $\theta$ acts on $K$. Let
$K_R$ be the fixed point set of $K$ under $\theta$. One can prove
that $K_R$ is the image under the projection to $P(E)$ of the set
of zeros of $Q_R$ on $E^\theta$.

 The conjugation $\theta$ also acts on the big cell
$K^\times=\{[y]\in K\;|\; y_{12}\neq 0\}$. According to
(\ref{pluckerbigcell}), this action is just the hermitian conjugate
$$\theta: A=\begin{pmatrix}\alpha&\beta
\\\gamma&\delta\end{pmatrix}\;\;\mapsto\,\;
A^\dag=\begin{pmatrix}\bar \alpha &\bar
\gamma
\\\bar\beta&\bar \delta\end{pmatrix}.$$
The fixed point set is $K_R^\times=\{A\in M_2(\C)\;|\; A^\dag=A\}$.

To the complex group $\rSO(6,\C)$ acting on $K$  there corresponds
the real form $\rSO(4,2)$ acting on $K_R$ and its spin group
$\rSU(2,2)$. The  hermitian form  on $\C^4$ left invariant by this
real form is
$$(u,v)=u^\dag F v,\qquad
F=\begin{pmatrix}0&I\\-I&0\end{pmatrix},\qquad u,v\in \C^4.$$ We can
now, as in the complex case, compute the subgroup of $\rSU(2,2)$
that leaves $K_R^\times$ invariant:
$$\begin{pmatrix}L&M\\NL&R\end{pmatrix}\begin{pmatrix}I\\
A\end{pmatrix}=\begin{pmatrix}L+MA\\NL+RA\end{pmatrix},
\qquad \begin{pmatrix}L&M\\NL&R\end{pmatrix}\in \rSU(2,2)$$ so
$L+MA$ must be invertible for all hermitian $A$. The result is that
the subgroup $P_R$ leaving $K_R^\times$ invariant is
$$P_R=\left\{\begin{pmatrix}L&0\\NL&{L^\dag}^{-1}\end{pmatrix},\:
N \hbox{ hermitian and }\det (L)\in \R
\right\},$$ and its action on $K_R^\times$ is
$$A\mapsto N+{L^\dag}^{-1}AL^{-1}.$$
We see that $P_R$ is the Poincar\'e group ($\det (L)=1$) times
dilations. This completes the interpretation of $K^\times_R$ as
the real Minkowski space.

\section{Homogeneous spaces for Lie supergroups \label{homogeneous}}

A Lie supergroup is a supermanifold having a group structure, that
is a multiplication map $\mu: G \times G \lra G$ and an inverse $i:
G \lra G$ satisfying the usual diagrams. As in the ordinary theory
to a Lie super group is associated a Lie superalgebra consisting of
the left invariant vector fields and identified with the tangent
space to the supergroup at the identity. For more details see
\cite{va} and the Appendix \ref{supergroups}.

We are interested in the construction of homogeneous spaces for
Lie supergroups. Let $G$ be a Lie supergroup and $H$ a closed Lie
subsupergroup. We want to define a supermanifold $X\equiv G/H$
with reduced manifold as $X_0=G_0/H_0$, and to construct a
morphism $\pi : G\longrightarrow X$ such that the following
properties are satisfied:

 \begin{enumerate}
 \item The reduction of $\pi$ is the natural map
$\pi_0: G_0\longrightarrow X_0$.

 \item $\pi$ is a submersion.
 \item There is an action $\beta$ from the left of $G$ on
$X$ reducing to the action of $G_0$ on $X_0$ and compatible with
the action of $G$ on itself from the left through $\pi$:

$$
\begin{CD}
G\times G@>\mu>>G\\@V1\times\pi VV @VV\pi V\\
G\times X@>\beta>>X\end{CD}
$$
 \item The pair $(X, \pi)$, subject to the properties 1,
2, and 3 is unique up to isomorphism. The isomorphism between two
choices is compatible with the actions, and it is also unique.
 \end{enumerate}

\subsection{The supermanifold structure on $X=G/H$}

Let $G$ be a Lie supergroup of dimension $m|n$ and $H$ a closed sub
supergroup of dimension $r|s$. Let $ \fg={\rm Lie}(G)$ and $\fh={\rm
Lie}(H)$. For each $Z\in \fg$, let $D_Z$ be the left invariant
vector field on $G$ defined by $Z$ (see the Appendix A for more
details). For $x_0\in G_0$ let $\ell_{x_0}$ and  $r_{x_0}$ be the
left and right translations of $G$ given by $x_0$. We denote by
$i_{x_0}=\ell_{x_0}\circ r_{x_0}^{-1}$ the inner automorphism
defined by $x_0$. It fixes the identity and induces the
transformation ${\rm Ad}_{x_0}$ on $\fg$.

For any open subset
 $U\subset G_0$ and any sub-superalgebra $\fk$ of $\fg$ we define
${\cO}_\fk(U)$ by
$$
{\oo}_\fk(U)=\{f\in {\oo}_G(U)\ |\ D_Zf=0 \hbox { on } U \hbox {
for all }Z\in \fk \}.
$$
Then ${\oo}_\fk$ is a subsheaf of ${\oo}_G$ (and, in particular,
so is ${\oo}_\fh$). On the other hand, for any open subset
$W\subset G_0$, invariant under right translations by elements of
$H_0$, we put
$$
{\oo}_{H_0}(W)=\{f\in {\oo}_G(W)\ |\ f \hbox { is invariant under
} r_{x_0} \hbox { for all } x_0\in H_0\}.
$$
 If $H_0$ is connected we have
$$
{\oo}_{H_0}(W)={\oo}_{\fh_0}(W).
$$
For any open set $W_0\subset X_0=G_0/H_0$ with $W=\pi_0^{-1}(W_0)$ we put
$$
{\oo}_X(W_0)={\oo}_{H_0}(W)\cap {\oo}_\fh(W).
$$
Clearly $ {\oo}_X(W_0)={\oo}_\fh(W) $ if $H_0$ is connected. The
subsheaf ${\oo}_X$ is a supersheaf on $X_0$. We write $X\equiv
G/H$ for the ringed superspace $(X_0, {\oo}_X)$. Our aim is to
prove that $X$ is a supermanifold and that it has all the
properties 1 to 4 listed previously.


It is clear that the left action of the group $G_0$ on $X_0$
leaves ${\oo}_X$ invariant and so it is enough to prove that there
is an open neighborhood $W_0$ of $\pi_0(1)\equiv\bar 1$ with the
property that $(W_0, {\oo}_X|_{W_0})$ is a super domain, i. e.,
isomorphic to an open submanifold of $k^{p|q}$.

We will do  this using the local Frobenius Theorem
\ref{frobenius}. Also, we identify $\fg$ with the space of all
left invariant vector fields on $G_0$, thereby identifying the
tangent space of $G$ at every point canonically with $\fg$ itself.

On $G_0$ we have a distribution spanned by the vector fields in
$\fh$. We denote it by ${\dd}_\fh$.

On each $H_0$-coset $x_0H_0$ we have a supermanifold structure
which is a closed sub supermanifold of $G$. It is an integral
manifold of ${\dd}_\fh$, i.~e. the tangent space at any point is
the subspace $\fh$ at that point. By the local Frobenius theorem
there is an open neighborhood $U$ of $1$ and coordinates $x_i$,
$1\leq i\leq n$ and $\theta _\alpha$, $1\leq \alpha\leq m$ on $U$
such that at each point of $U$, ${\dd}_\fh$ is spanned by
$\partial/\partial x_i ,
\partial/\partial \theta _\alpha (1\le i\le r, 1\le \alpha \le
s)$. Moreover, from the theory on $G_0$ we may assume that the
slices $L({\bf c}):=\{(x_1, \dots , x_n)\, |\;\; x_j=c_j,
\;\;r+1\le j\le n\}$ are open subsets of  distinct $H_0$-cosets
for distinct ${\bf c}=(c_{r+1},\dots ,c_n)$. These slices are
therefore supermanifolds with coordinates $x_i$, $\theta_\alpha$,
$1\le i\le r,\;\; 1\le \alpha \le s$. We have a  sub-supermanifold
$W'$ of $U$ defined by $x_i=0$ with  $1\le i\le r$ and
$\theta_\alpha =0$ with  $1\le \alpha \le s$. The map
$\pi_0:G_0\lra X_0$ may be assumed to be a diffeomorphism of
$W'_0$ with its image $W_0$ in $X_0$ and so we may view $W_0$ as a
superdomain, say $W$. The map $\pi_0$ is then a diffeomorphism of
$W'$ with $W$. What we want to show is that $W=(W_0,
{\oo}_X|_{W_0})$.

\begin{lemma} The map
$$
\begin{CD}
W'\times H@>\gamma>> G\\
 w,h@>>> wh\end{CD}
$$
is a super diffeomorphism of $W'\times H$ onto the open
sub-supermanifold of $G$ with reduced manifold the open subset
$W'_0H_0$ of $G_0$.
\end{lemma}

{\it Proof.} The map $\gamma $ in question is the informal
description of the map $\mu\circ (i_{W'}\times i_H)$ where $i_M$
refers to the canonical inclusion $M\hookrightarrow G$ of a
sub-supermanifold of $G$ into $G$, and $\mu:G\times
G\longrightarrow G$ is the multiplication morphism of the Lie
supergroup $G$. We shall use such informal descriptions without
comment from now on.

It is classical that the reduced map $\gamma_0$ is a
diffeomorphism of $W'_0\times H_0$ onto the open set $U=W_0'H_0$.
This uses the fact that the cosets $wH_0$ are distinct for
distinct $w\in W'_0$. It is thus enough to show that $d\gamma$ is
surjective at all points of $W'_0\times H_0$. For any $h\in H_0$,
right translation by $h$ (on the second factor in $W'\times H$ and
simply $r_h$ on $G$) is a super diffeomorphism commuting with
$\gamma$ and so it is enough to prove this at $(w,1)$.  If $X\in
\fg$ is tangent to $W'$ at $w$ and $Y\in \fh$, then
$$
d\gamma (X,Y)=d\gamma (X,0)+d\gamma (0,Y)=d\mu (X,0)+d\mu (0,Y)=X+Y.
$$
Hence the range of $d\gamma$ is all of $\fg$ since, from the
coordinate chart at $1$ we see that the tangent spaces to $W'$ and
$wH_0$ at $w$ are transversal and span the tangent space to $G$ at
$w$ which is $\fg$. This proves the lemma.\hfill$\blacksquare$

\begin{lemma}
We have
$$
\gamma^\ast {\oo}_X\big|_{W_0}={\oo}_{W'}\otimes 1,
$$
where $\gamma^\ast: \cO_G \lra \gamma^\ast \cO_{W' \times H}$.
\end{lemma}

{\it Proof.} To ease the notation we drop the open set in writing
a sheaf superalgebra, that is we will write $\cO_X$ instead of
$\cO_X(U)$.

We want to show that for any $g$ in ${\oo}_X\big|_U$, $\gamma^\ast
g$ is of the form $f\otimes 1$ and that the map $g\longmapsto f$
is bijective with ${\oo}_{W'}$. Now $\gamma^\ast$ intertwines $D_Z
(Z\in \fh)$ with $1\otimes D_Z$ and so $(1\otimes D_Z)\gamma^\ast
g=0$. Since the $D_Z$ span all the super vector fields on $H_0$ it
follows using charts that for any $p\in H_0$ we have $\gamma^\ast
g=f_p\otimes 1$ locally around $p$ for some $f_p\in {\oo}_{W'}$.
Clearly $f_p$ is locally constant in $p$. Hence $f_p$ is
independent of $p$ if $H_0$ is connected. If we do assume that
$H_0$ is connected, the right invariance under $H_0$ shows that
$f_p$ is independent of $p$. In the other direction it is obvious
that if we start with $f\otimes 1$ it is the image of an element
of ${\oo}_X\big|_U$.\hfill$\blacksquare$

\begin{theorem}
The superspace $(X_0, {\oo}_X)$ is a supermanifold.
\end{theorem}

{\it Proof}. At this stage by the previous lemmas we know that
$(X_0, {\oo}_X)$ is a super manifold at $\bar 1$. The left
invariance of the sheaf under $G_0$ shows this to be true at all
points of $X_0$. The proof that $(X_0, {\oo}_X)$ is a super
manifold is finished. \hfill$\blacksquare$

\subsection{The action of $G$ on $X$}

Clearly $G_0$ acts on $X$ but there is more: there is a natural
action of the supergroup $G$ itself on $X$. We shall now describe
how this action comes about.

\begin{proposition}
There is a map $ \beta  : G\times X\longrightarrow X $ such that
the following diagram
$$\begin{CD}G\times G@>\mu>>G\\@V1\times\pi VV @VV\pi V\\G\times X@>\beta>>X\end{CD}
$$
commutes.
\end{proposition}

{\it Proof.} Let $ \alpha :=\pi \circ \mu  : G\times
G\longrightarrow X$. The action of $G_0$ on $X_0$ shows that such
a map $\beta_0$ exists at the reduced level. So it is a question
of constructing the pull-back map
$$
\beta^\ast : {\oo}_X \longrightarrow {\oo}_{G\times X}
$$
such that
$$
(1\times \pi)^\ast \circ \beta^\ast=\alpha^\ast.
$$

Now $\pi^\ast$ is an {\it isomorphism\/} of ${\oo}_X$ onto the
sheaf ${\oo}_G$ restricted to a sheaf on $X$ ($W\longmapsto
{\oo}_G(\pi_0^{-1}(W))$), an so to prove the {\it existence and
uniqueness\/} of $\beta^\ast$ it is a question of proving that
$\alpha^\ast$ and $(1\times \pi)^\ast$ have the same image in
${\oo}_{G\times G}$. It is easy to see that $(1\times \pi)^\ast$
has as its image the subsheaf of sections $f$ killed by $1\otimes
D_X (X\in \fh)$ and invariant under $1\times r(h) (h\in H_0)$. It
is not difficult to see that this is also the image of
$\alpha^\ast$. \hfill$\blacksquare$

\medskip

We tackle now the  question of the uniqueness of $X$ (point 4 at
the beginning of Section \ref{homogeneous}).



\begin{proposition}
Let $X'$ be a super manifold with $X'_0=X_0$ and let $\pi'$ be a
morphism $G\longrightarrow X'$. Suppose that

\medskip

\noindent{(a)} $\pi'$ is a submersion.

\medskip

\noindent{(b)} The fibers of $\pi'$ are the super manifolds which
are the cosets of $H$.

\medskip

\noindent Then there is a natural isomorphism
$$
X\simeq X'.
$$
\end{proposition}

{\it Proof.} Indeed, from the local description of submersions as
projections it is clear that for any open $W_0\subset X_0$, the
elements of $\pi'^\ast ({\oo}_{X'}(W_0))$ are invariant under
$r(H_0)$ and killed by $D_X (X\in \fh)$. Hence we have a natural
map $X'\longrightarrow X$ commuting with $\pi$ and $\pi'$. This is
a submersion, and by dimension considerations it is clear that
this map is an isomorphism. \hfill$\blacksquare$

\medskip

We now turn to examine the functor of points of a quotient.

\subsection{The functor of points of $X=G/H$}

We are interested in a characterization of the functor of
points of a quotient $X=G/H$ in terms of the functor of points
of the supergroups $G$ and $H$.

\begin{theorem} \label{quotient}
Let $G$ and $H$ be supergroups as above and let $Q$ be the
sheafification 
of the functor: $T \lra G(T)/H(T)$. Then $Q$ is
representable and it is the functor of points of the homogeneous
space supermanifold $X=G/H$ constructed above.
\end{theorem}
{\it Proof.} In order to prove this result we can use the
uniqueness property which characterizes the homogeneous space $X$
described above. So we only need to prove the two facts: 1. $Q$ is
representable and 2. $\pi:G \lra Q$ is a submersion, the other
properties being clear.

To prove that $Q$ is representable we use the criterion in Theorem
\ref{representability}. The fact that $Q$ has the sheaf property
is clear by its very definition. So it is enough to prove there is
a open subfunctor around the origin (then by translation we can
move it everywhere). But this is given by $W \cong W' \times H$
constructed above. The fact $\pi$ is a submersion comes by looking
at it in the local coordinates given by $W$. \hfill$\blacksquare$

\begin{example} {\sl The flag supermanifold as a quotient.} \label{quotientflag}
Let $G=\rSL(m|n)$ the complex special linear supergroup as
described in Appendix \ref{supergroups}. There is a natural action
of $G$ on $V=\C^{m|n}$:
$$
\begin{CD}
G(T) \times V(T) @>>> V(T) \\
 (g_{ij}, v_k) @>>> \sum_j g_{ij}v_j
\end{CD}
$$

This action extends immediately to the flag supermanifold
$F=F(r|s,p|q; m|n)$ of $r|s$, $p|q$ subspaces in $\C^{m|n}$
described in detail in Example \ref{grassmannian}. Let us fix a
flag $\cF=\{\cO_T^{r|s} \subset \cO_T^{p|q}\}$ and consider the
map:
$$
\begin{CD}
G(T) @>>> F(T) \\
 g  @>>> g \cdot \cF.
\end{CD}
$$
The stabilizer subgroup functor is the subfunctor of $\rGL(4|1)$
given as
$$
H(T)=\{ g \in \rSL(m|n)(T) \;\; | \;\; g \cdot \cF = \cF\} \subset
\rSL(m|n)(T).
$$
Using the submersion Theorem one can see that this functor is
representable by a group supermanifold (see Section
\ref{superminkowski} for more details).

The flag supermanifold functor $F$ is isomorphic to the functor
$G/H$ defined as the sheafification of the functor $T \lra
G(T)/H(T)$. In fact we have that locally $G(T)$ acts transitively
on the set of direct summands of fixed rank $\cO_T$, hence the map
$g \lra g \cdot \cF$ is surjective. Hence by Theorem
\ref{quotient} we have that as supermanifolds $F \cong
G/H$.

The general flag is treated similarly. \hfill$\blacksquare$

\end{example}

\section{The super Poincar\'{e} and the translation superalgebras\label{spst}}

We want to give a description of the  Poincar\'{e}
superalgebra\footnote{For us the Poincar\'{e} algebra and Poincar\'{e}
superalgebra are meant to contain also the dilation generator.}
and the translation superalgebras, as subalgebras of the conformal
(Wess-Zumino) superalgebra. We will describe these superalgebras
in detail, over the complex and the real fields.

Results concerning these algebras have been known in the physics
literature for some time (see for example Ref. \cite{do}, the
reader can compare the approaches). In our exposition, we try to
extract the properties that uniquely select the Poincar\'{e}
superalgebra as the suitable supersymmetric extension of the
Poincar\'{e} algebra. This will be the natural setting to construct
the global theory in Section \ref{global}.

\subsection{The translation superalgebra}

We first describe the translation superalgebra as it comes from
physics.

The translation superalgebra is a  Lie superalgebra
$\fn=\fn_0+\fn_1$ such that $\fn_0$ is an abelian Lie algebra of
dimension $4$, which acts trivially on $\fn_1$, also of dimension
4. We have different descriptions depending if the ground field is
$\R$ or $\C$:

 \begin{enumerate}
 \item In the real case,  let $\fl_0$ be the Lie algebra of
the Lorentz group, $\fl_0=\fso(3,1)$. The subspace $\fn_1$ is a
real spin module for $\fl_0$. It is the sum of the two
inequivalent complex spin modules
 $\fn_1=S^+\oplus S^-$, and  satisfies a reality condition\footnote{On $S^+\oplus S^-$ there is a conjugation
  $\sigma$ commuting with the action of the spin group. This conjugation exchanges the spaces $S^+$ and $S^-$,
   and in a certain
  basis is just the complex conjugation.
  The spinors $s$ satisfying the reality condition $\sigma(s)=s$ provide an irreducible real representation
  of the spin group. They are called  Majorana spinors in the physics literature.}.
 \item In the complex case, $\fl_0$ is the complexification
of the Lorentz algebra, $\fl_0=\fsl_2\oplus\fsl_2$. The
subspace$\fn_1$ is, in physicists notation, the $\fl_0$-module
$\fn_1\simeq D(1/2,0)\oplus D(0,1/2)$.
 \end{enumerate}
 In both cases the odd commutator is a  non zero
symmetric $\fl_0$-map from $\fn_1\times \fn_1$ into $\fn_0$ which
is $\fl_0$-equivariant.

These properties identify uniquely the translation superalgebra in
both cases, real and complex.

The translation superalgebra will be identified in Section
\ref{superminkowski} with the dual of the infinitesimal Minkowski
superspace time. In Subsection \ref{complexpoincare} we will see
that it is the complement in the conformal superalgebra of the
Poincar\'{e} superalgebra.

\subsection{The complex Poincar\'{e} superalgebra\label{complexpoincare}.}

Let us restrict our attention to the complex field. The classical
complex Minkowski space time, as described in Section \ref{minkowski},
has a natural action of the Poincar\'{e} group naturally embedded into
the conformal group. Its compactification, $G(2,4)$, carries a
natural action of the complex conformal group $\rSL(4, \C)$.

 In
the super geometric infinitesimal setting the complex conformal
(Wess-Zumino) superalgebra is  $\fg=\fsl(4|1)$. We want to define
the Poincar\'{e} Lie superalgebra $\fp$ as a generalization of the Lie
algebra of the Poincar\'{e} group introduced in Section
\ref{minkowski}. We want $\fp$ to be a Lie subsuperalgebra in
$\fg$ subject to the following natural conditions:

 \begin{enumerate}
 \item $\fp$ is a  parabolic superalgebra, i.~e. it
contains a  Borel sub-superalgebra of $\fg$ (the Borel
subsuperalgebra is normalized so that its intersection with
$\fg_0$ is the standard Borel of $\fg_0$).
 \item $\fp\cap \fg_0$ is the parabolic subalgebra of
$\fg_0$ that consists of matrices of the form
$$
\begin{pmatrix} L&0&0\\ M&R&0\\0&0&0
\end{pmatrix},
$$
where $L, M, R\in \cM_2(\C)$ are $2\times 2$ matrices.
 \item There is a sub superalgebra $\fn \subset \fg$ such
that $\fg=\fp \oplus \fn$, and $\fn$ is a translation superalgebra
in the sense described above. The Lorentz algebra$\fl_0$ is the
subalgebra of $\fp\cap \fg_0$ consisting of matrices of the form
$$
\fl_0=\left\{\begin{pmatrix} L&0&0\\ 0&R&0\\0&0&0
\end{pmatrix}, \quad L, R\in \cM_2(\C)\right\},
$$
and the group $L_0$ is the simply connected group associated to
$\fl_0$, namely
$$
L_0=\left\{\begin{pmatrix} x&0&0\\ 0&y&0\\ 0&0&1\\
\end{pmatrix}, \quad x, y\in {\rSL}(2, {\C})\right\}.
$$
In particular this amounts to ask that $(\fg/\fp)_1$ as an
$L_0$-module is isomorphic to $D(1/2,0)\oplus D(0,1/2)$.
 \end{enumerate}
It turns out that these conditions determine uniquely the Poincar\'{e}
superalgebra $\fp$ inside the conformal superalgebra $\fg$.

\begin{lemma}\label{lemmapn}
There exists a unique subalgebra $\fp$ in the conformal
superalgebra $\fg$ satisfying the conditions 1 to 3 above, namely
$\fg=\fp \oplus \fn$ where
\begin{equation}
\fp=\left\{\begin{pmatrix} L&0&0\\ M&R&\alpha \\ \beta&0&c \\
\end{pmatrix}\right\}, \qquad
\fn=\left\{\begin{pmatrix} 0&A&\gamma\\ 0&0&0\\ 0&\delta&0\\
\end{pmatrix}\right\}.\label{superpoincare}
\end{equation}
Here $L, M, R$ and $A$ are $2\times 2$ matrices, $\gamma$ and
$\alpha $ are $1\times 2$, $\delta$ and $\beta$ are $2\times 1$,
and $c$ is a scalar and the supertrace condition for $\rSL(4|1)$
is $c={\rm tr}(L)+{\rm tr}(R)$.
\end{lemma}

{} From now on we refer to $\fp$ as the \textit{Poincar\'{e} Lie
superalgebra} or \textit{Poincar\'{e}  superalgebra} for shortness.

{\it Proof}. It is easy to check that $\fp$ and $\fn$ are Lie
superalgebras. This can be done by matrix calculation, but it is
easier (and better suited for the higher dimensional case) to look
at these as composed of root spaces with respect to the Cartan
subalgebra $\frak h$ where $\frak h$ consists of diagonal matrices
$$
H=\left\{\begin{pmatrix} a_1&0&0&0&0\\0&a_2&0&0&0\\ 0&0&a_3&0&0\\ 0&0&0&a_4&0\\ 0&0&0&0&a_5\\
\end{pmatrix}\right\}\qquad \hbox{with  }\; a_5=a_1+a_2+a_3+a_4.
$$
We treat the $a_i$ as the linear functions $H \longmapsto a_i$,
for $1\le i\le 5$. Then $\fp$ is the sum of $\fh$ and the root
spaces for the roots
\begin{eqnarray*}
&\pm (a_{1}-a_{2}),\; \pm (a_{3}-a_{4}),\;
a_3-a_i,\; a_4-a_i,\; a_5-a_i\;&\hbox{ for } i=1,2;\\
&a_j-a_5 &\hbox{ for } j=3,4; \end{eqnarray*} while $\fn_0, \fn_1$
are the respective sums of root spaces for the roots
$$
a_j-a_i,\; \hbox{ and } \; a_5-a_i,\;a_j-a_5\quad \hbox{ for }
i=3,4,\; j=1,2.
$$
The root description above implies easily that $0\not=[\fn_1,
\fn_1]\subset \fn_0$, and also that $\fn_0$ acts trivially on
$\fn_1$, i.e., $[\fn_0,\fn_1]=0$.

To verify the module structure of $\fn_1$ under $L_0$ it is more
convenient to use the matrix description. The formula
$$
\begin{pmatrix} x&0&0\\ 0&y&0 \\ 0&0&1\\ \end{pmatrix}
\begin{pmatrix} 0&0&\gamma\\0&0&0\\0&\delta&0\\ \end{pmatrix}
\begin{pmatrix} x^{-1}&0&0 \\ 0&y^{-1}&0 \\ 0&0&1\\ \end{pmatrix}=
\begin{pmatrix}
0&0&x\gamma \\ 0&0&0 \\ 0&\delta y^{-1}&0\\
\end{pmatrix}
$$
shows that the action is
\begin{equation}
(\gamma, \delta)\longmapsto (x\gamma,\delta
y^{-1})\label{actlorentz}
\end{equation}
which gives
$$
\fn_1\simeq D(1/2,0)\oplus D(0,1/2).
$$
There are $4$ other parabolic sub-superalgebras as one can see
examining the complete list in the general case in \cite{ka}
p.~51. They are defined as the sets of matrices of the following
four different forms:

\begin{eqnarray*}
&\fp_1=\left\{\begin{pmatrix}  a_{11} & a_{12} & 0 & 0 & 0  \\
                             a_{21} & a_{22} & 0 & 0 & 0 \\
                             a_{31} & a_{32}& a_{33} & a_{34} & 0 \\
                             a_{41} & a_{42} & a_{43} & a_{44} & 0\\
                             \beta_1 & \beta_2 & \beta_3& \beta_4 & a_{55} \\
\end{pmatrix}\right\},
&\fp_2=\left\{\begin{pmatrix}  a_{11} & a_{12} & 0 & 0 & 0  \\
                             a_{21} & a_{22} & 0 & 0 & 0 \\
                             a_{31} & a_{32}& a_{33} & a_{34} & 0 \\
                             a_{41} & a_{42} & a_{43} & a_{44} & \alpha_4\\
                             \beta_1 & \beta_2 & \beta_3& 0 & a_{55} \\
\end{pmatrix}\right\},\\
\\
&\fp_3=\left\{\begin{pmatrix}  a_{11} & a_{12} & 0 & 0 & 0  \\
                             a_{21} & a_{22} & 0 & 0 & \al_2 \\
                             a_{31} & a_{32}& a_{33} & a_{34} & \al_3 \\
                             a_{41} & a_{42} & a_{43} & a_{44} & \al_4\\
                             \beta_1 & 0 & 0& 0 & a_{55} \\
\end{pmatrix}\right\},
&\fp_4=\left\{\begin{pmatrix}  a_{11} & a_{12} & 0 & 0 & \al_1 \\
                             a_{21} & a_{22} & 0 & 0 & \al_2 \\
                             a_{31} & a_{32}& a_{33} & a_{34} & \al_3 \\
                             a_{41} & a_{42} & a_{43} & a_{44} & \al_4\\
                             0 & 0 & 0& 0 & a_{55} \\
\end{pmatrix}\right\}.
\end{eqnarray*}

However each they fail condition 3. The verifications are the same
as above and are omitted.

\begin{observation}
We can define a form on $\frak n_0$ by taking the determinant:
$$
q : \begin{pmatrix} 0&A&0\\ 0&0&0\\ 0&0&0\\
\end{pmatrix} \longmapsto \det (A).
$$
The odd commutator map is (using the notation of Lemma
\ref{lemmapn})
$$
\left((\gamma,\delta), (\gamma',\delta')\right)\longmapsto
\begin{pmatrix} 0&\gamma\delta'+\gamma'\delta&0\\ 0 &0&0\\ 0&0&0\\
\end{pmatrix},
$$
and under the transformation (\ref{actlorentz}) becomes
$$
A=\gamma\delta'+\gamma'\delta\longmapsto
x(\delta \gamma'+\delta'\gamma)y^{-1},
$$
showing that the form $q(\gamma,\delta)=\det (\gamma\delta'+\gamma'\delta)$
is preserved.\hfill$\blacksquare$
\end{observation}

\subsection{The real conformal and Poincar\'{e} superalgebras
\label{realsuperpoincare}}

We shall now show that this whole picture carries over to ${\R}$.
We shall construct a {\it conjugation\/} $\sigma$ of the super Lie
algebra $\fg=\fs\frak l(4|1)$ with the following properties:
\begin{enumerate}
 \item $\fg^\sigma$, the set of fixed points of $\sigma$,
is the Wess-Zumino super conformal Lie algebra, $\fsu(2,2|1)$.
 \item $\sigma$ leaves $\fp$ and $\fn$ invariant as well
as the even and odd parts of $\fg$.
 \item $\sigma$ preserves the big cell $\fc:=\fn\cap
\fg_0$.
 \item $\fc^\sigma$ consists of all $2\times 2$ skew
hermitian matrices.
 \item The group $L_0^\sigma$ consists of all matrices of
the form
$$
\begin{pmatrix} x&0\\ 0&{x^\dagger}^{-1}\\ \end{pmatrix},
$$
and the action on $\fc^\sigma$ is
$$
\begin{pmatrix} 0&A\\ 0&0\\ \end{pmatrix}\longmapsto
\begin{pmatrix} 0&xAx^\dagger\\ 0&0\\ \end{pmatrix}.
$$
\end{enumerate}
The function
$$
q : \begin{pmatrix}
0&A\\ 0&0\\
\end{pmatrix} \longmapsto \det (A)
$$
is the Minkowski metric; traditionally one takes the hermitian
matrices, but skew hermitian serve equally well: one just has to
multiply by $i$.

To construct $\sigma$ we proceed as in Ref. \cite{va}, pp.
112-113, but using
a slight variant of the construction. There, $F$ was taken in
block diagonal form as
$$
\begin{pmatrix} I_2&0\\ 0&-I_2\\ \end{pmatrix}.
$$
We now take
$$
F=\begin{pmatrix} 0&I_2\\ I_2&0\\ \end{pmatrix}
$$
where the blocks are all $2\times 2$ matrices.

\begin{proposition}
There is a conjugation $\sigma:\fg \lra \fg$ satisfying the
properties (1)-(5) listed above. It is given by
\begin{equation}
\sigma : \begin{pmatrix} X&\mu\\ \nu&x\\ \end{pmatrix} \longmapsto
\begin{pmatrix}
-FX^\dagger F&iF\nu^\dagger\\ i\mu^\dagger F&-\bar x\\
\end{pmatrix}\label{conjugation}
\end{equation}
where $X$ is  a $4\times 4$ matrix, $\mu$ is  $4\times 1$, $\nu$ is
$1\times 4$, and $x$ is a scalar.
\end{proposition}

{\it Proof.} It is a simple calculation to check that $ \sigma$ is
an antilinear map satisfying $\sigma^2=\id$ (so it is a
conjugation) and that it preserves the super bracket on $\fg$.  We
shall now verify Properties (1--5). In fact the verification of
Properties (2--5) is just routine and is omitted. The reason why
we take $F$ in the off diagonal form is to satisfy the requirement
that $\sigma$ leaves $\fp$ and $\fn$ invariant. That Property (1)
is also true is because the matrix $F$ defines a hermitian form of
signature $(2,2)$, as does the original choice. The fact that the
corresponding real form is transformable to the one where $F$ is
chosen in the diagonal form is proved on page 112 of Ref.
\cite{va}, with $F$ replaced by any hermitian matrix of signature
$(2,2)$. \hfill$\blacksquare$

\medskip

On the Poincar\'{e} superalgebra $\fp$ and on the super translation
algebra $\fn$ the conjugation $\sigma$ is given explicitly by:
\begin{equation}
\sigma :
\begin{pmatrix} L & 0 & 0 \\
                         M & R &\al \\
                         \be & 0 & c \\
\end{pmatrix} \longmapsto
\begin{pmatrix}
-R^\dagger & 0 & 0 \\
-M^\dagger & -L^\dagger & i\be^\dagger\\
i\al^\dagger & 0 &-\bar c\\
\end{pmatrix}\label{pconjugation}
\end{equation}

\begin{equation}
\sigma :
\begin{pmatrix}
0 & A & \ga \\
0 & 0 & 0 \\
0 & \de & 0 \\
\end{pmatrix} \longmapsto
\begin{pmatrix}
0 & A^\dagger & i\de^\dagger \\
0 & 0 & 0 \\
i\al^\dagger & 0 &-\bar c\\
\end{pmatrix}.\label{nconjugation}
\end{equation}

Hence the reality conditions read:
$$
R^\dagger=-L, \quad M=-M^\dagger, \quad \al=i\be^\dagger,
\quad c=-\bar c, \quad A= -A^\dagger, \quad \de=\ga^\dagger.
$$

\section{The global theory \label{global}}

We want to extend the infinitesimal results of the previous
section to obtain the Minkowski superspace time as the big cell
inside a certain super flag manifold, realized as homogeneous
space for the super conformal group.

\subsection{The complex super flag $F=F(2|0,2|1,4|1)$\label{superminkowski}}

Let $G=\rSL(4|1)$ be the complex super special linear group as
described in Appendix \ref{supergroups}. The natural action of $G$
on $V=\C^{4|1}$  extends immediately to the flag supermanifold (see
Example \ref{flag})
 $F=F(2|0,2|1; 4|1)$ of $2|0$, $2|1$ subspaces in
$\C^{4|1}$.
Let us fix a
flag $\cF=\{\cO_T^{2|0} \subset \cO_T^{2|1}\}$ and consider the
map:
$$
\begin{CD}
G(T) @>>> F(T) \\
 g  @>>> g \cdot \cF.
\end{CD}
$$
The stabilizer subgroup functor is the subfunctor of $G$ given as
$$
H(T)=\{ g \in G(T) \;\; | \;\; g \cdot \cF = \cF\} \subset G(T).
$$
One can readily check it consists of all matrices in $G(T)$ of the
form:
$$
\begin{pmatrix}
g_{11} & g_{12} & g_{13} & g_{14} & \gamma_{15} \\
g_{21} & g_{22} & g_{23} & g_{24} & \gamma_{25} \\
0 & 0 & g_{33} & g_{34} & 0 \\
0 & 0 & g_{43} & g_{44} & 0 \\
0 & 0 & \gamma_{53} & \gamma_{54} & g_{55} \\
\end{pmatrix}.
$$
This functor is representable by a group supermanifold by the
Submersion Theorem \ref{submersion}. In fact, as we shall see
presently, the map $g \mapsto g \cdot \cF$ is a submersion.


We wish to describe explicitly $F$ as the functor of points of
the supermanifold quotient of the supergroups $G$ and $H$ as it
was constructed in Section \ref{homogeneous}. The homogeneous space
$X=G/H$ constructed in Section \ref{homogeneous} is proven to be unique
once three conditions are verified, namely:

\begin{enumerate}
 \item The existence of a morphism $\pi:G \lra X$, such that
its reduction is the natural map $\pi_0: G_0\longrightarrow X_0$.
 \item $\pi$ is a submersion and the fiber of $\pi$ over
$\pi_0(1)$ is $H$. Since $\pi$ is a submersion  the fiber is well
defined as a super manifold.
 \item There is an action from the left of $G$ on $X$
reducing to the action of $G_0$ on $X_0$ and compatible with the
action of $G$ on itself from the left through $\pi$.

 \end{enumerate}
Conditions (1) and (3) are immediate in our case. The only thing
that we have to check is that the map
$$
\begin{CD}
G(T) @>\pi>> F(T)\\  g @>>> g \cdot \cF
\end{CD}
$$
is a submersion. We have to verify that at all topological points
$\pi$ has surjective differential. It is enough to do  the
calculation of the differential at the identity element $e \in
G(T)$.

 The calculation, being local, takes place inside the big
cell of $U\subset F$, which lies inside the product of big cells
$U_1 \times U_2$ of the Grassmannians $G_1=G(2|0;4|1)$,
$G_2=G(2|1;4|1)$.

We want to give local coordinates in $U_1 \times U_2$ and in $U$.
We know by Example \ref{grassmannian} that they are affine spaces.
So let us write, in the spirit of (2), local coordinates for $U_1$
and $U_2$ as
$$
\left(\begin{pmatrix}
I \\
A \\
\alpha
\end{pmatrix}, \;
\begin{pmatrix}
I & 0 \\
B & \beta \\
0 & 1
\end{pmatrix} \right)\in U_1(T) \times U_2(T), \qquad T \in
\smfld,
$$
where $I$ is the identity, $A$ and $B$ are $2 \times 2$ matrices
with even entries and $\al=(\al_1, \al_2)$, $\be^t=(\be_1,\be_2)$
are rows with odd entries.

An element of $U_1$  is inside $U_2$ if
\begin{equation}A=B+\beta\alpha,\label{twistor}\end{equation} so we can take as coordinates for a
flag in the big cell $U$ the triplet $(A,\alpha,\beta)$. We see
then  that $U$ is an affine $4|4$ superspace. Equation
(\ref{twistor}) is also known as {\it twistor relation}, see Ref.
\cite{ma1}.

 In these coordinates, $\cF$ becomes
$$
\left(\begin{pmatrix}
I \\
0 \\
0
\end{pmatrix}, \;
\begin{pmatrix}
I & 0 \\
 0& 0 \\
0 & 1
\end{pmatrix} \right)\approx (0,0,0).
$$

We want to write the map $\pi$ in these coordinates. In a suitable
open subset near the identity of the group we can take an element
$g\in G(T)$ as
$$g=\begin{pmatrix}g_{ij}&\gamma_{i5}\\\gamma_{5j}&g_{55}\end{pmatrix},
\qquad i,j=1,\dots 4.$$ Then, we can write an element $g \cdot \cF
\in G_1 \times G_2$ as:
\begin{equation}
\begin{pmatrix}
g_{11} & g_{12} \\
g_{21} & g_{22} \\
g_{31} & g_{32} \\
g_{41} & g_{42} \\
\gamma_{51} & \gamma_{52}
\end{pmatrix},
\begin{pmatrix}
g_{11} & g_{12} & \gamma_{15}\\
g_{21} & g_{22} & \gamma_{25} \\
g_{31} & g_{32} & \gamma_{35} \\
g_{41} & g_{42} & \gamma_{45}\\
\gamma_{51} & \gamma_{52} & g_{55}
\end{pmatrix}
\quad \approx \quad
\begin{pmatrix}
I \\
W Z^{-1} \\
\rho_1 Z^{-1}
\end{pmatrix},
\begin{pmatrix}
I & 0 \\
V Y^{-1}& (\tau_2-WZ^{-1}\tau_1) a \\
0 & 1
\end{pmatrix},\label{bigcell}
\end{equation}
where
\begin{eqnarray*}
&\rho_1=
\begin{pmatrix}
\gamma_{51} & \gamma_{52}
\end{pmatrix},\quad W=\begin{pmatrix}
g_{31} & g_{32} \\
g_{41} & g_{42}
\end{pmatrix},\quad Z=\begin{pmatrix}
g_{11} & g_{12} \\
g_{21} & g_{22}
\end{pmatrix}, \\&\tau_1=
\begin{pmatrix}
\gamma_{15} \\ \gamma_{25}
\end{pmatrix},\quad \tau_2=
\begin{pmatrix}
\gamma_{35} \\ \gamma_{45}
\end{pmatrix},\quad  d=(g_{55}-\nu Z^{-1}\mu_1)^{-1}\\&
V=W-g_{55}^{-1}\tau_2\rho_1,\quad Y=Z-g_{55}^{-1}\tau_1\rho_1.
\end{eqnarray*}
Finally the map $\pi$ in these coordinates  is given by:
$$
g \mapsto \left(W Z^{-1}, \rho_1 Z^{-1}, (\tau_2-WZ^{-1}\tau_1)
d\right).
$$
At this point one can compute the super Jacobian  and verify that
at the identity it is surjective. Let $\bar e$ be the image of the
identity in $F$, $\bar e\in U$.  It is easy to see that the
matrices of the form
$$\begin{pmatrix}0&0&0\\\tilde W&0&\tilde \tau_2\\\tilde
\rho_1&0&0\end{pmatrix}$$  map isomorphically onto the image of
$d\pi$. These are just the transpose of the matrices inside $\fn$
in (\ref{superpoincare}). We then obtain that the translation
superalgebra can be interpreted as the dual of the infinitesimal
Minkowski superspace.

It is not difficult to show that the subgroup of  $G=\rSL(4|1)$
that leaves the big cell invariant is the set of matrices in $G$
of the form \begin{equation}\begin{pmatrix}L & 0& 0\\
NL&R&R\chi\\d\varphi
&0&d\end{pmatrix},\label{supergrouppoincare}\end{equation} with
$L, N, R$ being $2\times 2$ even matrices, $\chi$ and odd $1\times
2$ matrix, $\varphi$ a $2\times 1$ odd matrix and $d$ a scalar.
This is then the complex Poincar\'{e} supergroup, whose Lie algebra is
$\fp$ in (\ref{superpoincare}). The action of the supergroup on
the big cell can be written as
\begin{eqnarray*}&&A\lra R(A+\chi\alpha)L^{-1}+N,\\&& \alpha\lra
d(\alpha +\varphi)L^{-1},\\&& \beta\lra
d^{-1}R(\beta+\chi).\end{eqnarray*} If the odd part is zero, then
the action reduces to the one of the Poincar\'{e} group on the Minkowski
space in (\ref{poincare}).

We then see that the big cell of the flag supermanifold
$F(2|0,4|0,4|1)$ can be interpreted as the complex super,
Minkowski space time, being the flag its {\it superconformal
compactification}.

\subsection{The real Minkowski superspace \label{realmink}}

We want to construct a real form of the Minkowski superspace that
has been described in Section \ref{superminkowski}.

We start by explicitly computing the real form of ${\rSL}(4|1)$ that
corresponds to the $\sigma$ defined in Section \ref{superpoincare}.
We shall compute it on $G\equiv {\rm GL}(4|1)$, but it is easy to
check that the conjugation so defined leaves $G_1\equiv {\rm
SL}(4|1)$ invariant.

In Appendix \ref{realform} it is explained that in order to obtain a
real form of $G$ we need a natural transformation $\rho$ from $G$ to
its complex conjugate $\bar G$.
Let $R$ be a complex ringed super space. We start by defining
$$
\begin{CD}{}
G(R) @>>> \bar G(R) \\ \\
g=\begin{pmatrix} D& \tau\\ \rho & d & \end{pmatrix}@>>>
g^\theta=\begin{pmatrix} D^\dagger&j \rho ^\dagger\\ j\tau^\dagger&\bar d\\
\end{pmatrix}
\end{CD}.
$$
Here we use $j$ for either $i$ or $-i$. To get the Lie algebra
conjugation $\sigma$ we shall eventually choose $j=i$, but at this
stage we need not specify which sign we take.

\begin{lemma} We have
$$
(h g)^\theta = g^\theta h^\theta.
$$
\end{lemma}

{\it Proof.}  It is important to note that we have taken the
following  convention: if $\theta$ and $\xi$ are odd variables,
then
\begin{equation}\overline{\theta\xi}=\bar\theta\bar\xi.
\label{convention}\end{equation} This
convention is opposed to the one used in physics, namely
$$\overline{\theta\xi}=\bar\xi\bar\theta,$$ but as it is explained in
Ref. \cite{dm}, (\ref{convention}) is the one that makes sense
functorially. According to this convention, then for matrices
$X,Y$ with {\it odd } entries
$$
(\overline {XY})^T=-(\bar Y)^T(\bar X)^T.
$$ Then the lemma results from direct
calculation.\hfill$\blacksquare$

\medskip

 We are ready to define the involution which gives the
real form of $G$. We will denote it by $\xi$:
$$
\begin{CD}
 G(R) @>\xi>> \bar G(R)  \\
g @>>> g^\xi & :=L(x^\theta)^{-1}L\end{CD} \qquad
L=\begin{pmatrix} F&0\\0&1
\end{pmatrix}, \quad F=\begin{pmatrix} 0&1\\1&0
\end{pmatrix}.
$$
We have that $(hg)^\xi=h^\xi g^\xi$ and that $(g^\xi)^\xi=g$, so it
is a conjugation.

\begin{lemma}
The map $x \mapsto x^\xi$ is a natural transformation. It defines a
ringed space involutive isomorphism $\rho: G \lra \bar G$ which is
$\C$-antilinear.
\end{lemma}

\textit{Proof}. This is a direct check.
\hfill$\blacksquare$

\medskip

\begin{proposition}
The topological space $G^\xi$ consisting of the points fixed by
$\rho$ has a real supermanifold structure and the supersheaf is
composed of those functions $f\in \cO_G$ such that $
\overline{\xi^*(f)}=f$.
\end{proposition}

\textit{Proof}. Immediate from the definitions in Appendix
\ref{realform}. \hfill$\blacksquare$

\medskip

We now wish to prove that the involution $\xi$ that gives the real
form for $G$ corresponds to the involution $\sigma$ constructed at
the Lie algebra level in Section \ref{realsuperpoincare}.

\begin{proposition} The conjugation
$\xi$ of $G$ defined above
induces on ${\rm Lie}(G)$ the  conjugation
$$\begin{CD}
\begin{pmatrix} X&\mu\\ \nu&x\\ \end{pmatrix} @>>>
\begin{pmatrix} -FX^\dagger F&-jF\nu^\dagger\\ -j\mu^\dagger F&-\bar x\\
\end{pmatrix}\end{CD},
$$ which for $j=-i$ coincides with $\sigma$ in
(\ref{conjugation}).
\end{proposition}

{\it Proof.} We have to compute the tangent map
at the identity, so we shall write $$\begin{pmatrix} D&\tau\\
\rho&d\\ \end{pmatrix}\approx \begin{pmatrix} \id&0\\
0&1\\ \end{pmatrix}+ \varepsilon\begin{pmatrix} X&\mu\\
\nu&x\\ \end{pmatrix},\qquad  \begin{pmatrix} X&\mu\\
\nu&x\\ \end{pmatrix}\in\fg,$$ up to first order in $\varepsilon$.
We now compute $x^\xi$ up to first order in $\varepsilon$,
$$x^\xi\approx\begin{pmatrix} \id&0\\
0&1\\ \end{pmatrix}+ \varepsilon L\begin{pmatrix} -X^\dagger&-j\nu^\dagger\\
-j\mu^\dagger&-\bar x\\ \end{pmatrix}L,$$ from which the result
follows.\hfill$\blacksquare$

We have defined a real form of $G=\rGL(4|1)$ which is also well
defined on $G_1=\rSL(4|1)$ and agrees with the real form of the
Lie superalgebra discussed in Section \ref{realsuperpoincare}.
Also, it is easy to check that it reduces to a conjugation on the
 Poincar\'{e} supergroup (\ref{supergrouppoincare}). We can compute
the conjugation on an element of the Poincar\'{e} supergroup:
\begin{eqnarray*}
&& g=\begin{pmatrix}L&0&0\\M&R&R\chi\\d\varphi&0&d\end{pmatrix},\quad
g^{-1}=\begin{pmatrix}
L^{-1}&0&0\\
-R^{-1}ML^{-1}+\chi\varphi L^{-1}&R^{-1}&-\chi d^{-1}\\-\varphi
L^{-1}&0&d^{-1}\end{pmatrix}\\\\
&&g^\xi=\begin{pmatrix}{R^\dag}^{-1}&0&0\\-{L^\dag}^{-1}M^\dag
{R^\dag}^{-1}-{L^\dag}^{-1}\varphi^\dag\chi^\dag&{L^\dag}^{-1}&-j{L^\dag}^{-1}\varphi^\dag\\-j{\bar
d}^{-1}\chi^\dag&0&{\bar d}^{-1}\end{pmatrix}.\end{eqnarray*} It
follows that the fixed points are those that satisfy the conditions:
\begin{equation} L={R^\dag}^{-1},\qquad \chi=-j\varphi^\dag,\qquad
ML^{-1}=-(ML^{-1})^\dag -j{L^\dag}^{-1}\varphi^\dag\varphi L^{-1}
,\label{realityconditions} \end{equation} which reduce to $\sigma$
(restricted to the super Poincar\'{e} subalgebra) in (\ref{conjugation})
at the Lie algebra level.

To get a more familiar form for the reality conditions, we observe
that the last equation in (\ref{realityconditions}) can be cast as
$$M'L^{-1}\equiv ML^{-1}+\frac 12 j{L^\dag}^{-1}\varphi^\dag\varphi
L^{-1},\qquad M'=-{M'}^\dag.$$ This is just an odd translation,
and amounts to multiply $g$ on the right by the group element
$$g'=\begin{pmatrix}\id&0&0\\-\frac 12 jR^{-1}{L^\dag}^{-1}\varphi^\dag\varphi L^{-1}&\id&0\\0&0&1\end{pmatrix}.$$

We want now to compute the real form of the big cell. The first
thing to observe is that the real form is well defined on the
quotient space $G/H$ (the superflag), where $H$ is the group
described in Section \ref{superminkowski} and consists of elements
in $G$ stabilizing a certain flag. In fact, one can check as we
did for the Poincar\'{e} supergroup that $H$ is stable under $\sigma$.

Notice that a
point of the big cell $(A,\alpha,\beta)$ can be represented by an
element of the group
$$g=\begin{pmatrix}\id&0&0\\A&\id&\beta\\\alpha&0&1\end{pmatrix},$$
since
$$g\begin{pmatrix}\id\\0\\0\end{pmatrix}=\begin{pmatrix}\id\\A\\\alpha\end{pmatrix},\qquad
g\begin{pmatrix}\id&0\\0&0\\0&1\end{pmatrix}=\begin{pmatrix}\id&0\\A&\beta\\\alpha&1\end{pmatrix}
\approx\begin{pmatrix}\id&0\\A-\beta\alpha&\beta\\0&1\end{pmatrix}.$$
 We first compute the inverse,
$$g^{-1}=\begin{pmatrix}\id&0&0\\-A+\beta\alpha&\id&-\beta\\-\alpha&0&1\end{pmatrix},$$
and then $g^\xi$
$$g^{\xi}=\begin{pmatrix}\id&0&0\\-A^\dag-\alpha^\dag\beta^\dag&\id&-j\alpha^\dag\\-j\beta^\dag&0&1\end{pmatrix}.$$
The element $g^\xi$ is already in the desired form, so the real
points are given by
$$A=-A^\dag-j\alpha^\dag\alpha,\qquad \beta=-j\alpha^\dag.$$
We can make a convenient change of coordinates,
$$A'\equiv A+\frac 12 j\alpha^\dag\alpha,$$ so the reality condition is
$$A'=-{A'}^\dag,$$ and we recover the same form than in Section (\ref{realsuperpoincare})
for the (purely even) Minkowski space time.

\appendix

\section{Supergeometry}\label{supergeometry}

In this section we want to recall the basic definitions and results
in supermanifold theory. For more details see Refs. \cite{ma1, dm,
va}, which use a language similar to ours.

\subsection{Basic definitions }\label{supergeometry1}

For definiteness, we take the ground field to be  $k=\R,\C$. A {\it
superalgebra} $A$ is a $\Z_2$-graded algebra, $A=A_0 \oplus A_1$.
The subspace $A_0$ is an algebra, while $A_1$ is an $A_0$-module.
Let $p(x)$ denote the parity of an homogeneous element $x$,
$$p(x)=0 \hbox{  if } x\in A_0,\qquad p(x)=1 \hbox{  if } x\in A_1.$$
The superalgebra $A$ is said to be {\it commutative} if for any two
homogeneous elements $x, y$
$$
xy=(-1)^{p(x)p(y)}yx
$$
The category of commutative superalgebras will be denoted by
$\salg$. From now on all superalgebras are assumed to be commutative
unless otherwise specified.

\begin{definition}
A {\it superspace} $S=(S_0, \cO_S)$ is a topological space $S_0$
endowed with a sheaf of superalgebras $\cO_S$ such that the stalk at
a point  $x\in S_0$ denoted by $\cO_{S,x}$ is a local
superalgebra\footnote{A local superalgebra is a superalgebra with a
maximal ideal.} for all $x \in S_0$. More generally we speak also of
\textit{ringed superspace} whenever we have a topological space and
a sheaf of superalgebras. \hfill$\blacksquare$
\end{definition}

\begin{definition} A {\it morphism} $\phi:S \lra T$ of superspaces is given by
$\phi=(\phi_0, \phi^*)$, where $\phi_0: S_0 \lra T_0$ is a map of
topological spaces and $\phi^*:\cO_T \lra \phi^*_0\cO_S$ is such
that $\phi_x^*(\bm_{\phi_0(x)})=\bm_x$ where $\bm_{\phi_0(x)}$ and
$\bm_{x}$ are the maximal ideals in the stalks $\cO_{T,\phi_0(x)}$
and $\cO_{S,x}$ respectively. \hfill$\blacksquare$
\end{definition}

\begin{example}
The superspace $k^{p|q}$ is the topological space $k^p$ endowed with
the following sheaf of superalgebras. For any open subset $U \subset
k^p$
$$
\cO_{k^{p|q}}(U)=\cO_{k^p}(U)\otimes k[\xi_1 \dots \xi_q],
$$
where $k[\xi_1 \dots \xi_q]$ is the exterior algebra (or {\it
Grassmann algebra}) generated by the $q$ variables $\xi_1 \dots
\xi_q$ and $\cO_{k^p}$ denotes the $C^{\infty}$ sheaf on $k^p$ when
$k=\R$ and the complex analytic sheaf on $k^p$ when $k=\C$ .
\hfill$\blacksquare$
\end{example}

\begin{definition}\label{supermanifold}
A  {\it supermanifold} of dimension $p|q$ is a superspace $M=(M_0,
\cO_M)$ which is locally isomorphic to $k^{p|q}$, i.~e. for all $x
\in M_0$ there exist an open set $V_x \subset M_0$ and $U \subset
k^{p}$ such that:
$$
{\cO_{M}}|_{V_x} \cong {\cO_{k^{p|q}}}|_U.
$$
A {\it morphism} of supermanifolds is simply a morphism of
superspaces.

The classical manifold $M_0$ underlying the supermanifold $M$, is
the {\it reduced space} of $M$. Its sheaf is given at any open set
$U$ as $\cO_M(U)$ modulo the nilpotent elements.
\hfill$\blacksquare$\end{definition}

The theory of supermanifolds resembles very closely the classical
theory. One can, for example, define tangent bundles, vector fields
and the differential of a morphism similar to the classical case.

\begin{definition}
Let $M=(M_0,\cO_M)$ be a supermanifold. A \textit{tangent vector}
$X_m$ at $m \in M_0$ is a (super) derivation $X_m:\cO_{M,m} \lra k$.
This allows usto define the tangent space $T_mM$ of $M$ at $m$ and
the tangent bundle $TM$. \textit{Super vector fields} are sections
of the tangent bundle. If $f:M \lra N$ is a morphism, we define its
differential $(df)_m:T_mM \lra T_{f(m)}N$ as $(df)_m (X_m) \alpha =
X_m(f^*(\alpha))$.

\hfill$\blacksquare$
\end{definition}

We summarize here some of the results that we will need later,
sending the reader to \cite{va}, Sections 4.2, 4.3, 4.4 and to
\cite{cf}, Chapter 4. These papers treat only the $C^\infty$ case;
the complex analytic case is done similarly.

\begin{definition}
Let $f:M \lra N$ be a supermanifold morphism with $f_0: M_0 \lra
N_0$ the underlying classical morphism on the reduced spaces. We say
that $f$ is a {\it submersion} at $m \in M_0$ if $f_0$ is an
ordinary submersion and $(df)_m$ is surjective. Likewise, $f$ is an
{\it immersion} at $m$ if $f_0$ is an immersion and $(df)_m$ is
injective. Finally $f$ is a {\it diffeomorphism} at $m$ if it is a
submersion and an immersion.

When we say $f$ is a submersion (resp. immersion or diffeomorphism)
we mean $f$ a submersion at all points of $M_0
$.\hfill$\blacksquare$
\end{definition}

As in the classical setting, submersions and immersions have the
usual local models. For more details see \cite{va} p.~148.

\begin{definition}
We say that $N$ is an {\it open } sub-supermanifold of $M$ if $N_0$
is an open submanifold of $M_0$ and $\cO_N$ is the restriction of
$\cO_M$ to $N$.

We say that $N$ is a {\it closed } sub-supermanifold of $M$ if $N_0$
is a closed submanifold of $M_0$ and there exists a map $f:N
\hookrightarrow M$ which is an immersion and such that $f_0$ is the
classical embedding $N_0 \subset M_0$ .\hfill$\blacksquare$
\end{definition}

Closed sub supermanifolds can be determined by using the super
version of the submersion theorem.

\begin{theorem}{\sl Submersion Theorem}.\label{submersion}
Let $f:M \lra N$ be a submersion and let $P_0=f_0^{-1}(n)$ for $n
\in N_0$. Then $P_0$ admits a supermanifold structure. Locally we
have that for $p \in P_0$, $\cO_{P, p}=\cO_{M,p}/f^*(I_n)$, where
$I_n$ is the ideal in $\cO_{N,n}$ of elements vanishing at $n$.
Moreover,
$$
\dim P= \dim M - \dim N.
$$
\end{theorem}

{\it Proof}. (Sketch). For an open subset $V$ of $P_0$, sections are
defined as assignments $q\longmapsto s(q)$ where $s(q)\in
{\oo}_{P,p}(q)$ for all $q$ and locally on $V$ these arise from
sections of ${\oo}_M$. By the local description of submersions as
projections it is seen easily that this defines the structure of a
super manifold on $P$. Note that $\dim (P)=\dim (M)-\dim (N)$.
\hfill$\blacksquare$

\medskip


Frobenius theorem plays a fundamental role in constructing sub
supermanifolds of a given manifold.

\begin{theorem}
Let $M$ be a supermanifold and let $\dd$ be an integrable super
distribution on $M$ of dimension $r|s$ (i.e. a subbundle, locally a
direct factor of the tangent bundle of $M$).

 {\sl Local Frobenius Theorem}. Then at each point there exists a
coordinate system $(x,\xi)$ such that the distribution at that point
is spanned by $\partial_{x_i}$, $\partial_{\xi_\alpha}$, with  $1
\leq i \leq r$ and  $1 \leq \alpha \leq s$.

 {\sl Global Frobenius Theorem}.\label{frobenius} Then at each
point there exists a unique maximal supermanifold $N$ such that
$TN=\dd$.
\end{theorem}

{\it Proof.} See Ref. \cite{va} 4.7, p.~157 and Ref. \cite{cf}
Chapter 4. \hfill$\blacksquare$

\medskip

We now turn to a different and alternative way to introduce
supermanifolds, that is the functor of points approach.

\subsection{The functor of points }\label{supergeometry2}

In supergeometry the functor of points approach is particularly
useful since it brings back the geometric intuition to the problems,
leaving the cumbersome sheaf notation in the background.

\begin{definition}
Given a supermanifold $X$ we define its \textit{functor of points}
as the following representable functor from the category of
supermanifolds to the category of sets:
$$
h_X:\smfld \lra \set, \qquad h_X(T)=\Hom(T,X).
$$
\end{definition}
\hfill$\blacksquare$

Given two supermanifolds $X$ and $Y$, Yoneda's lemma establishes a
one to one correspondence between the morphisms $X \lra Y$ and the
natural transformations $h_X \lra h_Y$. This allows one to view a
morphism of supermanifolds as a family of morphisms $h_X(T) \lra
h_Y(T)$ depending functorially on the supermanifold $T$.

We want to give a representability criterion, which allows one to
single out among all the functors from the category of
supermanifolds to the category of sets those that are representable,
i.~e.those that are the functor of points of a supermanifold. In
order to do this, we need to generalize the notion of {\it open
submanifold} and of {\it open cover} to fit this more general
functorial setting.

\begin{definition}. Let $U$ and $F$ be two functors $\smfld \lra \set$.
The functor $U$ is a \textit{subfunctor} of $F$ if $U(R) \subset
F(R)$ for all $R\in \smfld$. We denote it as $U\subset F$.

We say that $U$ is an \textit{open subfunctor} of $F$ if for all
natural transformations $f:h_T \lra F$ with $T \in \smfld$ then
$f^{-1}(U)=h_V$, where $V$ is open in $T$. If $U$ is also
representable we say that $U$ is an \textit{open supermanifold
subfunctor}.

Let $\cU_{\alpha}$ be open subfunctors of $\R^{m|n}$ (or
$\C^{m|n}$).
We say that $\{\cU_{\alpha}\}$ is an open cover of a functor
$F:\smfld \lra \set$ if for all supermanifolds $T$ and  natural
transformations $f:h_T \lra F$, $f^{-1}(\cU_\al)=h_{V_\al}$ and
$V_\al$ cover $T$.\hfill$\blacksquare$
\end{definition}

\begin{definition} A functor  $F:\smfld \lra \set$ is  said to be
{\it local} if it has the sheaf property, that is, if when
restricted to the open sets of a given supermanifold $T$ it is a
sheaf.
\end{definition}

Notice that any functor $F:\smfld \lra \set$ when restricted to the
category of open sub supermanifolds of a given supermanifold $T$
defines a presheaf.

We are ready to state a representability criterion which gives
necessary and sufficient conditions for a functor from $\smfld$ to
$\set$ to be representable.

\begin{theorem} {\sl Representability Criterion.} \label{representability}
Let $F$ be a functor $F: \smfld \lra \set$, such that when
restricted to the category of manifolds is representable.

\noindent  Then the functor $F$ is representable if and only if:

\noindent 1. $F$ is local, i.~e.it has the sheaf property.

\noindent 2. $F$ is covered by  open supermanifold functors.

\medskip

\end{theorem}

{\it Proof}. The proof resembles closely the  proof given in Ref.
\cite{dg} Chapter 1, for the ordinary algebraic category. For
completeness we include a sketch of it.

If $F$ is representable, $F=h_X$ and one can check directly that it
has the two properties listed above. This is done in the super
algebraic category for example in Ref. \cite{cf} Chapter 5.

Let $\{h_{X_\al}\}_{\al \in A}$ be the open supermanifold
subfunctors that cover $F$. Define $h_{X_{\al\be}}=h_{X_\al}
\times_F h_{X_\be}$ (This will correspond to the intersection of the
two open  $X_\al$ and $X_\be$).

We have the commutative diagram:
\begin{equation*}
\begin{CD}
h_{X_{\al\be}}=h_{X_\al} \times_F h_{X_\be}  @> j_{\be,\al} >>  h_{X_\be} \\
@VV{j_{\al,\be}}V @VV{i_\be}V\\
h_{X_\al} @> i_\al >>F
\end{CD}
\label{cd}
\end{equation*}

As a set we define:
$$
|X|=_{\defi} \coprod_{\al} |X_\al| / \sim,
$$
where $\sim$ is the following equivalence relation:
\begin{eqnarray*}
&&\forall\, x_\al \in |X_\al|, \,x_\be \in |X_\be|,\; x_\al \sim
x_{\be} \Longleftrightarrow \\&& \exists \, x_{\al\be} \in
|X_{\al\be}|,\,\hbox{ with }\,  j_{\al,\be}(x_{\al\be})=x_\al,\,
j_{\be,\al}(x_{\al\be})=x_{\be}
\end{eqnarray*}
and $|Y|$ denotes the underlying topological space of a generic
supermanifold $Y$.

This is an equivalence relation. The map $\pi_\al:|X_\al|
\hookrightarrow |X|$ is an injective map into the topological space
$|X|$.

We now need to define a sheaf of superalgebras $\cO_X$ by using the
sheaves in the open  $X_\al$ and ``gluing''. Let $U$ be open in
$|X|$ and let $U_\al=\pi_\al^{-1}(U)$. Define:
$$
\cO_X(U)=_{\defi}\{(f_\al) \in \coprod_{\alpha \in I}
\cO_{X_\al}(U_\al)| \quad
j_{\be,\gamma}^*(f_\be)=j_{\gamma,\be}^*(f_\gamma), \forall \be,
\gamma \in I \}.
$$
The condition $j_{\be,\gamma}^*(f_\be)=j_{\gamma,\be}^*(f_\gamma)$
simply states that to be an element of  $\cO_X(U)$, the collection
$\{f_\al\}$ must be such that $f_\be$ and $f_\gamma$ agree on the
intersection of $X_\be$ and $X_\gamma$ for any $\be$ and $\gamma$.

One can check directly that $\cO_X$ is a sheaf of superalgebras and
that $h_X=F$. For more details see \cite{cf} Chapter 5.
\hfill$\blacksquare$

\medskip

We want to discuss two important examples of supermanifolds together
with their functors of points, namely the  Grassmannian and the flag
supermanifolds. They will appear later in the definition of the
Minkowski super space time. In each case, we will define a functor
and then we will use Theorem \ref{representability} to show that the
functor is representable, so it is the functor of points of some
supermanifold whose reduced manifold is the Grassmannian and the
flag manifold respectively. We will see that our definitions of
super Grassmannian and super flag coincide with Manin's ones in Ref.
\cite{ma1}, Chapter 1.

\begin{remark}
For the rest of this paper we use the same letter, say $X$, to
denote both a supermanifold and its functor of points, as it is
customary to do in the literature. So instead of writing $h_X(T)$ we
will simply write $X(T)$, for $T$ a generic supermanifold.
\end{remark}

\begin{example}
 {\sl The Grassmannian supermanifold}.

\label{grassmannian} We define the Grassmannian of $r|s$-subspaces
of a $m|n$-dimensional complex vector space as the functor
$\rGr:\smfld \lra \set$ such that for any supermanifold $T$, with
reduced manifold $T_0$, $\rGr(T)$ is the set of locally free sheaves
over $T_0$ of rank $r|s$,  direct summands of $\cO_T^{m|n}=_{\defi}
\cO_T \otimes \C^{m|n}$.

Equivalently $\rGr(T)$ can also be defined as the set of pairs
$(L,\alpha)$ where $L$ is a locally free sheaf of rank $r|s$ and
$\alpha$ a surjective morphism
$$\alpha: \cO_T^{m|n} \lra L,$$ modulo the equivalence relation
$$(L,\alpha)\sim(L',\alpha')\quad \Leftrightarrow\quad L\approx
L', \quad \alpha'=a\circ\alpha,$$ where $a:L\rightarrow L$ is an
automorphism of $L$.

We need also to specify $\rGr$ on morphisms $\psi: R \lra T$.
Given the element $(L, \alpha)$ of  $\rGr(T)$,
$\alpha:\cO_T^{m|n}\rightarrow L$, we have the element of $\rGr(R)$,
$$
\rGr(\psi)(\alpha):\cO_R^{m|n}=\,\cO_T^{m|n}\otimes_{\cO_T}  \cO_R
\rightarrow L\otimes_{\cO_T} \cO_R.
$$

We want to show that $\rGr$ is the functor of points of a
supermanifold. By its very definition $\rGr$ is a local functor, so
by Theorem \ref{representability} we just have to show that $\rGr$
admits a cover by open supermanifold functors.

 Consider
the multiindex $I=(i_1, \dots, i_r| \mu_1, \dots, \mu_s)$ and the
map $\phi_I: \cO_T^{r|s} \lra \cO_T^{m|n}$ where
\begin{eqnarray*}\phi_I(x_1,\dots , x_r|\xi_1, \dots \xi_s)=&&m|n-\hbox{tuple
with }\\&&x_1, \dots , x_r \hbox{ occupying the position } i_1,
\dots, i_r,\\&&
 \xi_1, \dots , \xi_s \hbox{ occupying the position } \mu_1,
\dots, \mu_s \\&&\hbox{and the other positions are occupied by
zero}.\end{eqnarray*} For example, let $m=n=2$ and $r=s=1$. Then
$\phi_{1|2}(x,\xi)=(x,0|0, \xi)$.

We define the subfunctors $v_I$ of $\rGr$ as follows.  The set
$v_I(T)$ is the set of  pairs $(L,\al)$,  $\al: \cO_T^{m|n} \lra L$
(modulo the equivalence relation), such that $\al \circ \phi_I$ is
invertible. Since we can, up to an automorphism, choose
$\alpha(t)=t$ for any $t\in T_0$, this means that  $(\al \cdot
\phi)_t:\cO_{T,t}^{r|s} \lra L_t$ is an isomorphism of free
$\cO_{T,t}$-modules.

It is not difficult to check that they are open supermanifold
functors and that they cover $\rGr$. Actually the open supermanifold
functors $v_I$ are the functors of points of superspaces isomorphic
to  matrix superspaces of suitable dimension, as it happens in the
classical case.

This is very similar to the algebraic super geometry case which is
explained in detail in Ref. \cite{cf} Chapter 5. Hence, by Theorem
\ref{representability}, $\rGr$ is the functor of points of a
supermanifold, that we will call the \textit{ super Grassmannian of
$r|s$ subspaces into a $m|n$ dimensional space.}\hfill$\blacksquare$
\end{example}

\begin{example} {\sl The flag super manifold}.\label{flag}
Let $F: \smfld \lra \set$ be the functor such that for any
supermanifold $T$, $F(T)$ is the set of all flags $S_1\subset \cdots
\subset S_i\subset\cdots \subset S_k\subset \cO_T^{m|n}$, where the
$S_i$ are locally free sheaves of rank $d_i|e_i$, direct summands of
$\cO_T^{m|n}$. We want to show that $F$ is representable, i.~e.it is
the functor of points of what we call the super flag of
$d_1|e_1,\dots  d_k|e_k$ spaces into a $m|n$-dimensional complex
vector super space.

 The functor $F$ is clearly local, so by Theorem \ref{representability} we just have  to show that
it admits a cover by open submanifold functors. For concreteness, we
will consider the super flag of  $d_1|e_1, d_2|e_2$ spaces into and
$m|n$ space, the general case being a simple extension of this.
Consider the natural transformation $\phi:F \lra \rGr_1 \times
\rGr_2$ given by
$$
\begin{CD}
\phi_T: F(T) @>>> \rGr_1(T) \times \rGr_2(T) \\
S_1 \subset S_2 @>>> (S_1, S_2),
\end{CD}
$$
where the $\rGr_i$ are functors of points of super Grassmannians.

Let $v_I^i$ be the open subfunctors of $\rGr_i$ described in Example
\ref{grassmannian}. Consider the subfunctors $u_{IJ}=\phi^{-1}(v_I^1
\times v_J^2)$. These are open. In fact if we have a natural
transformation $\psi:h_X \lra F$, then
$$\psi^{-1}(u_{IJ})=(\phi \circ \psi)^{-1}(v_I^1 \times v_J^2)=h_W$$ for $W$ open
in $X$ (the product of open sets is open in our topology). The fact
that they cover $F$ is also clear. As for the representability,
since the $v_I$ are the functors of points of super spaces
$\C^{p|q}$ for suitable $p$ and $q$, it is not hard to see that the
$u_{IJ}$ will be given by algebraic relations in the $v_I$'s
coordinates, hence they are also representable. By Theorem
\ref{representability} we have that $F$ is also representable.

\hfill$\blacksquare$
\end{example}

\subsection{The linear supergroups \label{supergroups}}

For completeness, we include here a brief summary of  the
definition of a Lie supergroup, some of its properties,  and the
construction of the supergroups $\rGL(m|n)$ and $\rSL(m|n)$. For
more details see \cite{va} and \cite{cf}.

The ground field is always  $k=\R$ or $\C$.

\begin{definition}  A \textit{Lie supergroup} is a supermanifold
whose functor of points
$$
 G :\smfld \lra \set
$$
is group valued.\hfill$\blacksquare$
\end{definition}

\begin{remark}
Saying that $G$ is group valued is equivalent to have the
following natural transformations:

\begin{enumerate}
 \item Multiplication $\mu:  G  \times  G  \lra  G $, such
that $\mu \circ (\mu \times \id)=(\mu \times \id) \circ \mu$,
i.~e.
$$
\begin{CD}
G \times G \times G @> \mu \times \id >> G \times G \\
@V{\id \times \mu}VV @VV \mu V \\
  G \times G @> \mu >> G  \\
\end{CD}
$$

 \item Unit $e:e_k \lra  G $, where $e_k:\smfld
\lra \set$, such that $\mu \circ (\id \otimes e)= \mu \circ (e
\times \id)$, i. e.
$$
\begin{CD}
G \times e_k @>{\id \times e} >>
G \times G @>{e \times \id} >> e_k \times G \\
         @VVV  @V{\mu}VV @VVV \\
             G@=    G@= G  \\
\end{CD}$$
\item Inverse $i: G  \lra  G $,
such that $\mu \circ (\id \times i)=e \circ \id$, i. e.
$$
\begin{CD}
G @> (\id, i)>> G \times G \\
@V{}VV @VV \mu V \\
  e_k @>e>> G  \\
\end{CD}
$$
 \end{enumerate}

\end{remark}

\begin{example}
Let $k^{m|n}=(k^m, \cO_{k^{m|n}})$ denote the supermanifold whose
reduced space is the affine space $k^m$ and with supersheaf:
$$
\cO_{k^{m|n}}(U)=\cO_{k^m}(U) \otimes \wedge(\xi_1, \dots, \xi_n),
\qquad U \subset k^m
$$
introduced in Appendix \ref{supergeometry}. Its functor of points is
given by:
$$
\smfld \lra \set, \qquad k^{m|n}(T)=\Hom(T, k^{m|n}).
$$
The set $k^{m|n}(T)$ can be identified with the set of
$m|n$-tuples $(t_1, \dots, t_m|\theta_1, \dots ,\theta_n)$, where
the $t_i$'s and $\theta_j$'s are, respectively, even and odd
global sections of $\cO_{T}$.

The set $k^{m|n}(T)$ is an additive group for all $T$, so
$k^{m|n}$ is a Lie supergroup.\hfill$\blacksquare$
\end{example}

Given a Lie supergroup one defines its Lie superalgebra
$\mathrm{Lie}(G)$ as the set of left invariant vector fields,
together with a natural super bracket on them.

We recall that a left invariant vector field is defined as
a vector field $V$ satisfying the condition:
$$
(V \otimes \id)i^*=i^*V
$$
where $\mu$ denotes the multiplication and $i(x,y)=\mu(y,x)$.

As in the ordinary case $\mathrm{Lie}(G)$ can be identified with the
tangent space of $G$ at the identity (which is a topological point).
For more details see \cite{cf} Chapter 4, and \cite{va} p.~276.

\begin{example}

Let  $\rGL(m|n): \smfld \lra \set$ be the functor such that
$\rGL(m|n)(T)$ are the invertible $m|n \times m|n$ matrices with
entries in $\cO_T(T)$:
\begin{equation}\begin{pmatrix}p_{m\times m}&q_{m\times n }\\r_{n\times
m}&s_{n\times n}
\end{pmatrix},\label{morphisms}\end{equation}
where the submatrices $p$ and $s$ have even entries and $q$ and
$r$ have odd entries. The invertibility condition implies that $p$
and $q$ are ordinary invertible matrices.

$\rGL(m|n)$ is the functor of points of a supermanifold (i.~e. it
is representable) whose reduced space is an open set $U$ in the
ordinary space $k^{{m^2}+n^2}$, namely the matrices with
invertible diagonal blocks. In fact one can readily check that the
supersheaf $\cO_{\rGL(m|n)}$ is
$$
\cO_{\rGL(m|n)}=\cO_{k^{m^2+n^2|2mn}}|_U
$$
Notice: the supermanifold $k^{m^2+n^2|2mn}$ can be identified with
the supermatrices $m|n \times m|n$.

The subfunctor $\rSL(m|n)$ of $\rGL(m|n)$ consists on all matrices
with Berezinian \cite{be} equal to $1$, where
\begin{equation}
\mathrm{Ber} \begin{pmatrix}p_{m\times m}&q_{m\times n
}\\r_{n\times m}&s_{n\times n}
\end{pmatrix}=\det(s^{-1})\det(p-qsr).
\end{equation}

The proof that $\rSL(m|n)$ is representable uses Theorem
\ref{submersion} (the submersion theorem).

It is not hard to show that $\rLie(\rGL(m|n))(T)$ consists of all
matrices $m|n \times m|n$ (with entries in $\cO_T(T)$) while
$\rLie(\rSL(m|n))(T)$ is the subalgebra of $\rLie(\rGL(m|n))$
consisting of matrices with zero supertrace.
\end{example}

\subsection{Real structures and real forms \label{realform}}

We want to understand how it is possible to define real structures
and real forms in supergeometry. For more details see Ref.
\cite{dm} p.~92.

A major character in this game is the complex conjugate of
a super manifold.

Let  $M=(|M|, \cO_M)$ be a complex manifold. The \textit{complex
conjugate} of $M$ is the manifold $\bM=(|M|, \cO_\bM)$ where
$\cO_\bM$ is the sheaf of the antiholomorphic functions on $M$
(which are immediately defined once we have $\cO_M$ and the complex
structure on $M$). We have a $\C$-antilinear sheaf morphism
$$\begin{CD}\cO_M@>>> \cO_\bM\\ f@>>> \bar f.\end{CD}$$

In the super context it is not possible to speak directly of antiholomorphic
functions and for this reason we need the following generalization of
complex conjugate super manifold.

\begin{definition}
Let $M=(|M|, \cO_M)$ be a complex super manifold. We define a
\textit{complex conjugate} of $M$ as a complex super manifold
$\bM=(|\bM|, \cO_{\bM})$, where now $\cO_\bM$ is just a
supersheaf, together with a ringed space $\C$-antilinear
isomorphism. This means that we have an isomorphism of topological
spaces $|M| \cong |\bM|$ and a $\C$-antilinear sheaf isomorphism
$$\begin{CD}\cO_M@>>> \cO_\bM\\ f@>>> \bar
f.\end{CD}$$\hfill$\blacksquare$
\end{definition}


\begin{example}
Let $M=\C^{1|1}=(\C, \cO_\C[\theta])$,  $\bM=(\C, \cO_{\bar\C}[\bar \theta])$
where $\cO_\C$ and $\cO_{\bar \C}$ denote
respectively the sheaf of holomorphic and antiholomorphic
functions on $\C$. The isomorphism is
$$
\begin{array}{ccc}
\cO_M & \lra & \cO_\bM \\
z & \mapsto & \bz \\
\theta & \mapsto & \btheta.
\end{array}
$$
Notice that while $\bar z$ has a the meaning of being the complex
conjugate of $z$, $\bar \theta$ is simply another odd variable
that we introduce to define the complex conjugate.

Practically one can think of the complex conjugate super manifold as
a way of giving a meaning to $\bar f$ the complex conjugate of a
super holomorphic function.\hfill$\blacksquare$
\end{example}


We are ready to define a real structure on a complex
supermanifold.

\begin{definition}
Let $M=(|M|, \cO_M)$ be a complex super manifold. We define a \textit{real
  structure} on $M$ as an involutive isomorphism of ringed spaces
$\rho: M \lra  \bM$, which is $\C$-antilinear on the sheaves
$\rho^*:\cO_\bM \lra \rho^* \cO_M$, $\rho^*(\lambda f)=\bar \lambda
\rho(f)$. We define the \textit{real form} $M_r$ of $M$ defined by
$\rho$ as the supermanifold $(M^\rho, \cO_{M_r})$ where $M^\rho$ are
the fixed points of $\rho:|M| \lra |\bM|=|M|$ and $\cO_{M_r}$ are
all the functions $f \in \cO_M|_{M^\rho}$ such that $
\overline{\rho^*(f)}=f$.\hfill$\blacksquare$
\end{definition}

If $M$ is a complex supermanifold, one can always construct the
complex conjugate $\bar M$ in the following way. Take $|\bM|=|M|$
and as $\cO_\bM$ the sheaf with the complex conjugate $\C$-algebra
structure (that is $\lambda \cdot f=\bar \lambda f$). In order to
obtain a real structure on $M$, we need a ringed spaces morphism
$M \lra \bM$ with certain properties. By Yoneda's Lemma this is
equivalent to give an invertible  natural transformation between
the functors of points:
$$
\rho:M(R) \lra \bM(R)
$$
for all super ringed spaces $R$ satisfying the $\C$-antilinear
condition.

We take this point of view in Section \ref{realmink} when we discuss
the real Minkowski space.

\section*{Acknowledgments}

We would like to thank Dr. Claudio Carmeli of INFN, Genoa, Italy,
for pointing out an imprecision in our original characterization of
the super homogeneous space in Section \ref{homogeneous}, as well
for helpful discussions.

 R. Fioresi thanks  Prof. A. Vistoli for helpful comments.

R. Fioresi and V. S. Varadarajan  want to thank the Universidad de
Valencia, Departamento de F\'{\i}sica Te\'{o}rica, for its kind hospitality
during the realization of this work.

M. A. Lled\'{o} wants to thank the Universit\`{a} di Bologna for its kind
hospitality during the realization of this work.

R. Fioresi wants to thank the University of California at Los
Angeles for its kind hospitality during the realization of this
work.

 This work has been supported by the the Spanish Ministerio de
Educaci\'{o}n y Ciencia through the grant FIS2005-02761 and EU FEDER
funds, by the Generalitat Valenciana, GV05/102 and and by the EU
network MRTN-CT-2004-005104 `Constituents, Fundamental Forces and
Symmetries of the Universe'.


\begin{thebibliography}{99}

\bibitem{ma1} Y. Manin. {\it Gauge field theory and complex
geometry.} Springer Verlag, (1988). (Original Russian edition in
1984).

\bibitem{pe} R. Penrose. {\it  Twistor algebra}. J. Math. Phys. {\bf 8}, 345-366, (1967).


\bibitem{be} F. A. Berezin. {\it
Introduction to superanalysis.} Edited by A. A. Kirillov. D. Reidel
Publishing Company, Dordrecht (Holland) (1987).  With an Appendix by
V. I. Ogievetsky. Translated from the Russian by J. Niederle and R.
Koteck\'y. Translation edited by Dimitri Le\u\i tes.


\bibitem{cf} L. Caston and  R. Fioresi. {\it  Mathematical Foundation of
Supersymmetry}. To appear.


\bibitem{dflv} R. D'Auria, S. Ferrara, M. A. Lled\'{o} and V. S.
Varadarajan.
{\it  Spinor algebras}. J.Geom.Phys. {\bf 40}, 101-128, (2001).

\bibitem{dm} P.  Deligne and J. Morgan.
{\it Notes on  supersymmetry (following J. Bernstein)}, in ``Quantum
fields and strings. A
            course for mathematicians", Vol 1, AMS, (1999).

\bibitem{dg} M. Demazure and P. Gabriel. {\it Groupes
Alg\'ebriques, Tome 1.} Mason $\&$ Cie, \'editeur. North-Holland
Publishing Company, The Netherlands (1970).

\bibitem{fe} A. Ferber.
{\it  Supertwistors and conformal supersymmetry.} Nucl. Phys. B,
{\bf 132}, 55-64, (1978).

\bibitem{ha} R. Hartshorne {\it  Algebraic geometry}. (Third edition). Springer
Verlag. New York, (1983).

\bibitem{ka} V. G. Kac. {\it Lie superalgebras}. Adv. in Math. {\bf
26},  8-26, (1977).

\bibitem{fwz} S. Ferrara, J. Wess and B. Zumino. {\it Supergauge
multiplets and superfields}. Phys. Lett. B {\bf 51}, 239, (1974).

\bibitem{kn} M. Kotrla and  J. Niederle. {\it  Supertwistors And Superspace.} Czech. J.Phys. B {\bf 35}, 602,
(1985).



\bibitem{hh} P. S. Howe and G. G. Hartwell. {\it  A superspace survey}.
Class. Quantum Grav. {\bf 12}, 1823-1880, (1995).


\bibitem{va} V. S. Varadarajan. {\it Supersymmetry for mathematicians: an
introduction}. Courant Lecture Notes, {\bf 1}.  AMS (2004).

\bibitem{do} V.K. Dobrev and V.B. Petkova. {\it On the group-theoretical
approach to extended conformal supersymmetry : function space
realizations and invariant differential operators}.  Fortschr. d.
Phys. {\bf 35}, 537-572, (1987).







%
%
%
%
%
%
%
%
%
%
%

%
%
%








\end{thebibliography}
\end{document}